\documentclass[twoside]{mystyle}

\newtheorem{theorem}{Theorem}[subsection]
\newtheorem{corollary}[theorem]{Corollary}
\newtheorem{lemma}[theorem]{Lemma}

\newcommand{\qed}{\mbox{$\Box$}}
\renewcommand{\qed}{\hskip0.25em\raisebox{0.6ex}{\framebox[0.5em][l]{\ }}\vspace{.5pc}}
\newcommand{\onespace}{\renewcommand{\baselinestretch}{1}\normalsize}
\newcommand{\twospace}{\renewcommand{\baselinestretch}{1.5}\normalsize}

\newcommand{\num}{\renewcommand{\theequation}{\thetheorem}\addtocounter{theorem}{1}}

\def\theequation{\thetheorem}

\NONUMBIB

\begin{document}

\twospace

\SPECFNSYMBOL{}{}{}{}{}{}{}{}{}%

\AOPMAKETITLE

\AOPAMS{Primary 60B15, 60J20; secondary 20E22.}
\AOPKeywords{Bernoulli--Laplace diffusion, Markov chain, hyperoctahedral group, homogeneous space, Fourier transform.}
\AOPtitle{A SIGNED GENERALIZATION OF THE BERNOULLI--LAPLACE DIFFUSION MODEL}
\AOPauthor{Clyde~H.~Schoolfield,~Jr.}
\AOPaffil{Harvard University}
\AOPlrh{C.H.~SCHOOLFIELD,~JR.}
\AOPrrh{SIGNED GENERALIZATION OF BERNOULLI--LAPLACE DIFFUSION MODEL}
\AOPAbstract{We bound the rate of convergence to stationarity for a signed generalization of the Bernoulli--Laplace diffusion model; 
this signed generalization is a Markov chain on the homogeneous space ($\mathbb{Z}_2~\wr~S_n) / (S_r~\times~S_{n-r})$.
Specifically, for $r$ not too far from $n/2$, we determine that, to first order in $n$, $\frac{1}{4} n \log n$ steps are both
necessary and sufficient for total variation distance to become small.  Moreover, for $r$ not too far from $n/2$, we show that our
signed generalization also exhibits the ``cutoff phenomenon.''}

\maketitle

\BACKTONORMALFOOTNOTE{3}

\thispagestyle{empty}

\twospace

\section{Introduction.} \label{1}

Consider the classical Bernoulli--Laplace model for the diffusion of gases through a membrane, in which at each step two randomly 
chosen balls from different urns are switched.  How many steps does it take for this process to achieve near-randomness?  This 
question was answered by Diaconis and Shahshahani (1987).  Suppose that the balls also have charges and that, at each step, the 
two balls are not only switched, but their charges are also possibly changed.  How many steps does it take for this process to 
achieve near-randomness?  This is the question that we consider.

For the Bernoulli--Laplace model, Diaconis and Shahshahani (1987) obtained bounds on the rate of convergence to
stationarity.  Similarly, in preparation for our main result, we bound the rate of convergence to stationarity for a variant of the
classical Bernoulli--Laplace diffusion model; this variant is also a Markov chain on the homogeneous space $S_n / (S_r \times
S_{n-r})$.  Specifically, for $r$ not too far from $n/2$, we determine that, to first order in $n$, $\frac{1}{4} n \log n$ steps are
both necessary and sufficient for total variation distance to become small.

We then bound the rate of convergence to stationarity for a signed generalization of our variant of the Bernoulli--Laplace
diffusion model; this generalization is a Markov chain on the homogeneous space ($\mathbb{Z}_2~\wr~S_n) / (S_r~\times~S_{n-r})$.
Specifically, for $r$ not too far from $n/2$, we determine that, to first order in $n$, $\frac{1}{4} n \log n$ steps are both
necessary and sufficient for total variation distance to become small.  Moreover, for $r$ not too far from $n/2$, we show that our
signed model also exhibits the ``cutoff phenomenon.''  We also examine a slight variant of this signed model, establishing upper and
lower bounds on its rate of convergence to stationarity.  

\section{The Bernoulli--Laplace Diffusion Model.} \label{2}

\subsection{Introduction.} \label{2.1}

We now review the Bernoulli--Laplace model for the diffusion of gases through a membrane.  This is done in preparation for
Section~\ref{3} where we extend the results of this section to a signed generalization of the Bernoulli--Laplace model.
Imagine two racks, the first with positions labeled $1$ through $r$ and the second with positions labeled $r+1$ through $n$.  Without
loss of generality, we assume that   $1 \leq r \leq n/2$.  Suppose that there are $n$ balls, labeled $1$ through $n$, each
initially placed at its corresponding position.  

At each step, a ball is chosen at random from each rack and the two balls are switched.  Then, if necessary, the balls on each of
the two racks are permuted so that their labels are in increasing order.  The state of the system is completely determined, at each
step, by the labels of balls on the first rack.  This switching procedure thus induces a Markov chain on the set of all ${n \choose
r}$ subsets of size $r$ from a set with $n$ elements.  

Let $K$ be the subgroup of $S_n$ which permutes the first $r$ indices among themselves and the last $n-r$ indices among themselves,
but does not commingle these two sets of indices.  Thus $K \cong S_r~\times~S_{n-r}$.  Notice that $K$ is the stabilizer of the
initial configuration of the process.  The switching procedure thus also induces a Markov chain on the homogeneous space $X = S_n /
(S_r~\times~S_{n-r})$.  The set $X$ may thus be identified with the set of all ${n \choose r}$ subsets of size $r$ from a set with
$n$ elements.

Let $T$ be the set of all transpositions in $S_n$.  Also let $T_1$ be the set of all transpositions in $K$ and let $T_2$ be the set
of all transpositions in $S_n \setminus K$.  Thus $T = T_1 \cup T_2$.  At each step, the process defined above chooses a random
element of $T_2$.

But notice that before the two balls to be switched have been chosen, the balls in the first rack may be permuted among themselves
and the balls in the second rack may be permuted among themselves, without changing the state of the system following the switch. 
Thus, at each step, the procedure actually chooses a random element of $T_2 K := \{ \tau \kappa \in S_n  : \tau \in T_2$ and $\kappa
\in K \}$.  Notice that each element of $T_2 K$ can be uniquely written as $\tau \kappa$, where $\tau \in T_2$ and $\kappa \in K$.

This Bernoulli--Laplace diffusion model may be modeled formally by a probability measure $P$ on the symmetric group $S_n$.  Since
$|T_2 K| = |T_2| \cdot |K| = r(n-r) \cdot r!(n-r)!$, we may thus define the following probability measure on the set of all
permutations of $S_n$:
\num \begin{equation} \label{2.1.1}
\begin{array}{rcll}

P(\tau\kappa) & := & \displaystyle \frac{1}{r(n-r)r!(n-r)!} & \mbox{where $\tau \kappa \in T_2 K$ and}
\vspace{.5pc} \\

P(\pi) & := & 0 & \mbox{otherwise}.

\end{array}
\end{equation}

\noindent
Since there are $n!$ elements in $S_n$, the uniform probability measure on the set of all permutations of $S_n$ is given by
\num \begin{equation} \label{2.1.2}
U(\pi) := \frac{1}{n!} \ \ \ \mbox{for every $\pi \in S_n$}.
\end{equation}

The following result, which is Theorem 2 of Diaconis and Shahshahani (1987), establishes an upper bound on both the total variation
distance and the $\ell^2$ distance between $\widetilde{P^{*k}}$ and $\widetilde{U}$, where $\widetilde{P^{*k}}$ is the probability
measure on the homogeneous space $X = S_n / (S_r~\times~S_{n-r})$, induced by the convolution $P^{*k}$ of $P$ with itself $k$
times, and $\widetilde{U}$ is the uniform probability measure on $X$.  (Homogeneous spaces and probability measures induced on them
are discussed in Section~\ref{2.2}.)  It should be noted that $\| \widetilde{P^{*k}} - \widetilde{U} \|_{\mbox{\rm \scriptsize TV}}$ 
is the total variation distance to uniformity after $k$ steps for the Markov chain on $X$ defined above, started at the chosen initial
configuration.  

\begin{theorem} \label{2.1.3}
Let $P$ and $U$ be the probability measures on the symmetric group $S_n$ defined in (\ref{2.1.1}) and (\ref{2.1.2}),
respectively.  Let $\widetilde{P^{*k}}$ be the probability measure on the homogeneous space $X = S_n / (S_r~\times~S_{n-r})$ induced
by $P^{*k}$ and let $\widetilde{U}$ be the uniform probability measure on $X$.  Let $k = \frac{1}{2} r \left( 1 - \frac{r}{n}
\right) \left( \log n + c \right)$.  Then there exists a universal constant $a > 0$ such that

\[ \| \widetilde{P^{*k}} - \widetilde{U} \|_{\mbox{\rm \scriptsize TV}} \ \ \leq \ \ \mbox{$\frac{1}{2} {n \choose r}^{1/2}$} \|
\widetilde{P^{*k}} - \widetilde{U} \|_2 \ \ \leq \ \ ae^{-c/2} \ \ \ \mathrm{for\ all\ } \mbox{$c > 0$}. \]

\end{theorem}

In the following sections we present the results that were needed to prove this theorem and which are used to prove analogous
results in Section~\ref{2.5} and in Section~\ref{3}.  In Section~\ref{2.2} we present the definitions and results necessary to study
Markov chains on homogeneous spaces.  In Section~\ref{2.3} we decompose the representation $L(X)$ (defined in Section~\ref{2.2}) of 
$S_n$ into its irreducible constituents.  In Section~\ref{2.4} we calculate the Fourier transform of the probability measure defined 
in (\ref{2.1.1}), using a procedure different from that of Diaconis and Shahshahani (1987), and show how it was used to prove 
Theorem~\ref{2.1.3}.  In Section~\ref{2.5} we perform a similar analysis on a variant of the classical model; this variant will be 
generalized in Section~\ref{3}. 

\subsection {Markov Chains on Homogeneous Spaces.} \label{2.2}

We now present basic properties and results regarding homogeneous spaces and Markov chains on them.  A more detailed introduction to
this subject may be found in Chapter 3 of Diaconis (1988).

An \emph{action} of a group $G$ on a set $X$ is a mapping from the Cartesian product $G~\times~X$ to $X$, with the image of $(g,
x)$ under this action being denoted by $gx$, which satisfies (i) $ex = x$ for the identity element $e \in G$ and all $x \in X$, and 
(ii) $(g_2 g_1)x = g_2(g_1 x)$ for all $g_1, g_2 \in G$ and  $x \in X$.  For an element $x \in X$, the set of elements $\{ g \in G : 
gx = x\}$ is called the \emph{stabilizer} of $x$; the stabilizer of an element $x \in X$ is a subgroup of $G$.  A group $G$ acts 
\emph{transitively} on a set $X$ if, for every $x_1, x_2 \in X$, there exists some $g \in G$ such that $gx_1 = x_2$.  A set with a
group acting transitively on it is called a \emph{homogeneous space}.

Suppose that $G$ acts transitively on a set $X$ and that $K$ is the stabilizer of some fixed element $x_0 \in X$.  The group $G$
acts on the left coset space $G/K$ by sending $(g, hK)$ to $(gh)K$ for all $g \in G$ and $hK \in G/K$.  The set $X$ and the left
coset space $G/K$ are isomorphic under this action.  We may thus identify $X$ with a set of left coset representatives
$\{x_0, x_1, \ldots, x_{m-1}\}$ of $K$ in $G$, where $x_0 = e \in G$ and $x_i \in G$ for $1 \leq i \leq m-1$.

A probability measure $P$ defined on a finite group $G$ induces a probability measure $\widetilde{P}$ on the set $X = G/K$ by
defining $\widetilde{P}(x_i) := P(x_iK)$ for $0 \leq i \leq m-1$, where $x_i K := \{ x_i k \in G  : k \in K \}$.  The transition
matrix \textbf{$\widetilde{\mbox{P}}$} of the Markov chain induced by the probability measure $\widetilde{P}$ is defined, for the
element at the intersection of the row corresponding to $x_i \in X$ and the column corresponding to $x_j \in X$, by

\[ \widetilde{\mbox{\textbf{P}}}_{x_i, x_j} \ \ := \ \ P(x_j K x_i^{-1}) \]

\noindent
where $x_j K x_i^{-1} := \{ x_j k x_i^{-1} \in G  : k \in K \}$.  According to Lemma 3 in Section F of Chapter 3 of Diaconis
(1988), \textbf{$\widetilde{\mbox{P}}$} is doubly stochastic.  Furthermore, the $k$-step transition probabilities 
\textbf{$\widetilde{\mbox{P}}^k$} for the Markov chain satisfy

\[ \widetilde{\mbox{\textbf{P}}}^k_{x_i, x_j} \ \ = \ \ P^{*k}(x_j K x_i^{-1}) \]

\noindent
and therefore, since $x_0 := e$,

\[ \widetilde{\mbox{\textbf{P}}}^k_{x_0, x_j} \ \ = \ \ P^{*k}(x_j K) \ \ = \ \ \widetilde{P^{*k}}(x_j). \]

\noindent
This confirms the statement concerning $\| \widetilde{P^{*k}} - \widetilde{U} \|_{\mbox{\rm \scriptsize TV}}$ that immediately precedes
Theorem~\ref{2.1.3}.

Let $V$ be a finite-dimensional vector space over the complex numbers and let \textbf{GL}$(V)$ be the \emph{general linear group}
of isomorphisms of $V$ onto itself.  Let $\rho: G \longrightarrow\ $\textbf{GL}$(V)$ be an \emph{irreducible representation} of $G$ with 
\emph{character} $\chi_{\rho}$ and \emph{dimension} $d_{\rho}$.  Suppose that $H$ is a subgroup of $G$ and that $\rho_1$ and $\rho_2$
are representations of $G$ and $H$, respectively.  Let $\rho_1 \downarrow_H^G$ be the representation of $H$ known as the \emph{restricted
representation} and let $\rho_2 \uparrow_H^G$ be representation of $G$ known as the \emph{induced representation}.

Let $L(X)$ be the set of all (complex-valued) functionals defined on $X$.  The group $G$ acts on $L(X)$ via the definition $(gf)(x)
:= f(g^{-1}x)$ for $g \in G$, $x \in X$, and $f \in L(X)$.  Since for fixed $g \in G$ this action is a bijective linear mapping of
$L(X)$ onto itself, $L(X)$ can also be regarded as a representation of $G$.  In fact, as a representation, $L(X)$ is isomorphic
to $\rho_0 \uparrow_K^G$, the \emph{trivial representation} $\rho_0$ of $K$ induced up to a representation of $G$.  According to Theorem 2 
in Section 1.4 of Serre (1977), every representation of a group $G$ is the direct sum of irreducible representations of $G$.  Thus $L(X)$
is the direct sum of irreducible representations of $G$.  

According to Theorem 6 in Section 2.5 of Serre (1977), the characters $\chi_1, \chi_2, \ldots, \chi_s$ of the irreducible
representations of a group $G$ form an orthonormal basis for the Hilbert space of class functions on $G$ with respect to the inner
product defined by

\[ \left\langle \psi, \chi \right\rangle_G := \frac{1}{|G|} \sum_{g \in G} \psi(g) \overline{\chi(g)}. \]

\noindent
For any irreducible representation $\rho$ of $G$, define

\[ m_{\rho} \ \ := \ \ \left\langle \chi_{\rho}, \chi_{L(X)} \right\rangle_G. \] 

\noindent
Thus $m_{\rho}$ is the multiplicity of the irreducible representation $\rho$ of $G$ in the decomposition of $L(X)$ into its
irreducible constituents.

The following useful result relating induced and restricted characters is the \emph{Frobenius reciprocity formula}, which is
Theorem 13 in Section 7.2 of Serre (1977).

\begin{lemma} \label{2.2.1}
Suppose that $H$ is a subgroup of $G$ and that $\psi$ and $\chi$ are characters of $H$ and $G$, respectively.  Then

\[ \left\langle \psi \uparrow_H^G , \chi \right\rangle_G \ \ = \ \ \left\langle \psi , \chi \downarrow_H^G \right\rangle_H, \]

\noindent
where the inner product on the left (resp., right) is calculated in $G$ (resp., $H$).
\end{lemma}

\noindent
It thus follows from Lemma~\ref{2.2.1} that 

\[ m_{\rho} \ \ = \ \ \left\langle \chi_{\rho}, \chi_{L(X)} \right\rangle_G \ \ = \ \ \left\langle \chi_{\rho}, \chi_{\rho_0}
\uparrow_K^G \right\rangle_G \ \ = \ \ \left\langle \chi_{\rho} \downarrow_K^G, \chi_{\rho_0} \right\rangle_K. \]

\noindent
Thus $m_{\rho}$ is also the multiplicity of the trivial representation $\rho_0$ of $K$ in the decomposition of the restriction 
of $\rho$ to $K$.  Furthermore, we have the following, which is Lemma 5 in Section F of Chapter 3 of Diaconis (1988).

\begin{lemma} \label{2.2.2}
The multiplicity $m_{\rho}$ of an irreducible representation $\rho: G \longrightarrow\ $ \emph{\textbf{GL}}$(V)$ of $G$ in the
decomposition of $L(X)$ into its irreducible constituents equals the dimension of the space of $K$-fixed vectors in $V$, i.e.,

\[ m_{\rho} \ \ = \ \ \mbox{\rm dim}\{ v \in V  : \rho(k)v = v \ \mathrm{for\ all}\mbox{$ \ k \in K \}$}. \]

\end{lemma}

For notational purposes, define $\widehat{I}(\rho) \ := \ I \oplus 0$, where $I$ is the $m_{\rho}$-dimensional identity matrix and
$0$ is the $(d_{\rho} - m_{\rho})~\times~(d_{\rho} - m_{\rho})$ zero matrix.  Thus $\widehat{I}(\rho)$ is the $d_{\rho} \times
d_{\rho}$ matrix \onespace $\left[ \begin{array}{cc} I & 0 \\ 0 & 0 \\ \end{array} \right]$.  \twospace Notice that tr$\left(
\widehat{I}(\rho) \widehat{I}(\rho)^* \right) = m_{\rho}$.  The preceding lemma leads to the following result, which will be useful
in the calculation of the Fourier transform in Sections~\ref{2.4} and \ref{3.3}.

\begin{lemma} \label{2.2.3}
Suppose that a finite group $G$ acts transitively on a finite set $X$ and that $K$ is the stabilizer of some fixed element $x_0 \in
X$.  Also suppose that $\rho: G \longrightarrow\ $ \emph{\textbf{GL}}$(V)$ is an irreducible representation of $G$.  Then there is
an orthonormal basis of $V$ such that

\[ \frac{1}{|K|} \sum_{\kappa \in K} \rho(\kappa) \ \ = \ \ \widehat{I}(\rho). \]

\end{lemma}

\proof{Proof} Lemma~\ref{2.2.2} asserts that $m_{\rho}$ is the dimension of the space of $K$-fixed vectors in $V$.  Choose an 
orthonormal basis in $V$ such that the first $m_{\rho}$ basis vectors are $K$-fixed.  It then follows from Theorem 1 in Section 1.3 of
Serre (1977) that, under this choice of basis, the representation $\rho$ splits as $V = V_1 \oplus V_2$, where dim$\ V_1 = m_{\rho}$
and dim$\ V_2 = d_{\rho} - m_{\rho}$. Thus $\rho \downarrow_K^G = \rho_1 \oplus \rho_2$, where $\rho_1$, as a representation on $V_1$,
is the direct sum of $m_{\rho}$ trivial representations on $K$ and $\rho_2$ is a representation on $V_2$. 

Notice that $\displaystyle \frac{1}{|K|} \sum_{k \in K} \rho(k)$ is the Fourier transform of the uniform distribution on $K$ at the
representation $\rho \downarrow_K^G$ of $K$.  Furthermore, for any finite group, the Fourier transform of any probability distribution
at the trivial representation is $1$ and the Fourier transform of the uniform distribution at any nontrivial representation is the
zero matrix.  It is from this that the desired result follows.  \qed  \vspace{.5pc}

When $m_{\rho} \leq 1$ for every irreducible representation $\rho$ of $G$, the decomposition of $L(X)$ is called \emph{multiplicity
free}.  In such a case the pair $(G, K)$ is called a \emph{Gelfand pair}.  Thus if $(G, K)$ is a Gelfand pair and $\rho$ is an
irreducible representation of $G$ occurring in $L(X)$, then $\widehat{I}(\rho)$ is the $d_{\rho}~\times~d_{\rho}$ matrix whose only
nonzero entry is a ``1'' in the $(1,1)$ position.  A more detailed introduction to Gelfand pairs may be found in Diaconis (1988) or
Macdonald (1995).

A probability measure $P$ defined on a finite group $G$ is called \emph{right K-invariant} if $P(gk) = P(g)$ for all $g \in G$ and
$k \in K$.  There is an analogous definition for left $K$-invariance.  A probability measure $P$ is called \emph{bi-K-invariant} if
it is both left and right $K$-invariant, i.e., if $P(k_1 g k_2) = P(g)$ for all $g \in G$ and all $k_1, k_2 \in K$.  Our probability
measure $P$ defined in (\ref{2.1.1}) is clearly right $K$-invariant.  Notice for any $\tau \kappa \in T_2 K$ and any $\kappa_1 \in
K$ that $\kappa_1 \tau \kappa = (\kappa_1 \tau \kappa_1^{-1}) (\kappa_1 \kappa)$ with $\kappa_1 \tau \kappa_1^{-1} \in T_2$ and
$\kappa_1 \kappa \in K$.  It follows that $P$ is also left $K$-invariant and hence bi-$K$-invariant.

According to Theorem 9 in Section F of Diaconis (1988), if $(G, K)$ is a Gelfand pair and $P$ is a bi-$K$-invariant probability
measure defined on $G$, then the Fourier transform $\widehat{P}(\rho)$ at any irreducible representation $\rho$ of $G$ is a constant
times $\widehat{I}(\rho)$, where $\widehat{I}(\rho)$ is the $d_{\rho}~\times~d_{\rho}$ matrix whose only nonzero entry is a ``1'' in
the $(1,1)$ position. 

We are now able to bound the distance to uniformity of a probability measure $\widetilde{P}$ induced on a homogeneous space $X =
G/K$ from a right $K$-invariant probability measure $P$ defined on a finite group $G$ in terms of the Fourier transform of $P$ by
use of the following, which is the Upper Bound Lemma in Section F of Chapter 3 in Diaconis (1988).

\begin{lemma} \label{2.2.4}
Suppose that a finite group $G$ acts transitively on a finite set $X$ and that $K$ is the stabilizer of some fixed element $x_0 \in
X$.  Also suppose that $P$ is a right $K$-invariant probability measure defined on $G$, that $\widetilde{P}$ is the induced
probability measure on the set $X$, and that $\widetilde{U}$ is the uniform probability measure on $X$.  Then

\[ \| \widetilde{P} - \widetilde{U} \|_{\mbox{\rm \scriptsize TV}}^2 \ \ \leq \ \ \mbox{$\frac{1}{4}$} |X| \cdot \| \widetilde{P} -
\widetilde{U} \|_2^2 \ \ = \ \ \mbox{$\frac{1}{4}$} \sum_{\rho} d_{\rho}\
\mathrm{tr}\mbox{$\left(\widehat{P}(\rho)\widehat{P}(\rho)^*\right)$} \]

\noindent
where the sum is taken over all \emph{nontrivial} irreducible representations of $G$ that occur at least once in $L(X)$.
\end{lemma}

\noindent
Notice in the special case $K = \{e\}$ that $X = G/K = G$ and that the preceding lemma reduces to the Upper Bound Lemma for groups
(see, e.g., Diaconis (1988), Chapter 3, Section B). 

\subsection{Irreducible Representations in $L(X)$.} \label{2.3}

As a representation of $S_n$, $L(X)$ is isomorphic to $\rho_0 \uparrow_{S_r~\times~S_{n-r}}^{S_n}$, where $\rho_0$ is the trivial
representation of $S_r~\times~S_{n-r}$.  The decomposition of $L(X)$ into its irreducible constituents is given by the following,
which is Lemma 2.2.19 of James and Kerber (1981).

\begin{lemma} \label{2.3.1}
Let $L(X)$ be the set of all (complex-valued) functionals defined on $X = S_n / (S_r~\times~S_{n-r})$.  Then, as a
representation of $S_n$,

\[ \displaystyle L(X) \ \ \cong \ \ \rho_{[n]} \ \oplus \ \rho_{[n-1,1]} \ \oplus \ \cdots \ \oplus \ \rho_{[n-r,r]} \]

\noindent
where $\rho_{[n-i,i]}$ is the irreducible representation of $S_n$ corresponding to the partition $[n-i,i]$ of $n$, for $0 \leq i
\leq r$.  Furthermore, the dimension of $\rho_{[n-i,i]}$ is given by

\[ \displaystyle d_{[n-i,i]} \ \ = \ \ {n \choose i} \ - \ {n \choose i-1} \]

\noindent
with $[n-0,0] \equiv [n]$ and the the usual conventions ${k \choose 0} = 1$ and ${k \choose -1} = 0$.
\end{lemma}

Notice that the decomposition of $L(X)$ is multiplicity free.  Thus $(S_n, S_r~\times~S_{n-r})$ is a Gelfand pair.  So for any
irreducible representation $\rho$ occurring in the decomposition of $L(X)$, $\widehat{I}(\rho)$ is the $d_{\rho}~\times~d_{\rho}$
matrix whose only nonzero entry is a ``1'' in the $(1,1)$ position.

\subsection{Analysis of the Classical Model.} \label{2.4}

In order to apply the Upper Bound Lemma (\ref{2.2.4}) to our bi-$K$-invariant probability measure $P$, we must now calculate the
Fourier transform at each nontrivial irreducible representation of $S_n$ occurring in the decomposition of $L(X)$.  Our calculations
are done with the aid of the following two lemmas.

In the special case when $P$ is a class function, it is a consequence of Schur's Lemma that the Fourier transform may
be calculated easily by use of the following, which is Lemma 5 of Diaconis and Shahshahani (1981).

\begin{lemma} \label{2.4.1}
Suppose that $\rho$ is an irreducible representation of a finite group $G$ with character $\chi$ and that $P$ is a class function.
For each conjugacy class $i$, let $P_i$ be the constant value of $P$ on the class, let $n_i$ be the cardinality of the class, and
let $\chi_i$ be the constant value of $\chi$ on the class.  Then the Fourier transform of $P$ is given by

\[ \widehat{P}(\rho) \ \ = \ \ \left[ \frac{1}{d_{\rho}} \sum_{i=1}^s P_i n_i \chi_i \right] I, \]

\noindent
where $d_\rho$ is the dimension of $\rho$, $I$ is the $d_\rho$-dimensional identity matrix, and the sum is taken over distinct
conjugacy classes.  
\end{lemma}

The following formulas, found in Section D of Chapter 3 and Section B of Chapter~7, respectively, of Diaconis (1988) are used to
calculate the numerical value of the Fourier transform.

\begin{lemma} \label{2.4.2}
Suppose that $\rho$ is an irreducible representation of $S_n$ corresponding to the partition $[\lambda] = [\lambda_1, \ldots,
\lambda_k]$ of $n$.  Let $r(\lambda) := \chi_{[\lambda]}(\tau) / d_{[\lambda]}$ with $\tau \in S_n$.  Then

\[ r(\lambda) \ = \ \displaystyle \frac{1}{n(n-1)} \sum_{j=1}^k \left[ \lambda_j^2 \ - \ (2j-1) \lambda_j \right] \ \ \
\mbox{$\mathrm{and}$} \ \ \
d_{[\lambda]} \ = \ \displaystyle n! \ \mbox{$\mathrm{det}$} \left( \frac{1}{(\lambda_i - i + j)!} \right)_{1 \leq i,j \leq k}, \]

\noindent
with $1/m! := 0$ if $m < 0$.

\end{lemma}

We now calculate the Fourier transform at each nontrivial irreducible representation of $S_n$ occurring in the decomposition of
$L(X)$.  Diaconis and Shahshahani (1987) did this with the aid of spherical functions.  Our technique is different and was used by
Greenhalgh (1989) and Scarabotti (1997) in their analyses of variants of the Bernoulli--Laplace diffusion model.

\begin{lemma} \label{2.4.3}
Let $P$ be the probability measure on $S_n$ defined in (\ref{2.1.1}).  Let $\rho = \rho_{[n-i,i]}$ (for some $1 \leq i
\leq r$) be a nontrivial irreducible representation of $S_n$ occurring in the decomposition of $L(X)$.  Then, in a certain basis, 
the Fourier transform is

\[ \widehat{P}(\rho) \ \ = \ \ \left[ 1 \ - \ \frac{i(n-i+1)}{r(n-r)} \right] \widehat{I}(\rho) \]

\noindent
where $\widehat{I}(\rho)$ is the $d_{\rho}~\times~d_{\rho}$ matrix whose only nonzero entry is a ``\ $1$'' in the $(1,1)$ position.
\end{lemma}

\proof{Proof} Recall that $\rho_{[n]}$ is the trivial representation of $S_n$.  Thus we must calculate the Fourier transform for the
other $r$ nontrivial irreducible representations of $S_n$ occurring in the decomposition of $L(X)$, which were found in 
Lemma~\ref{2.3.1}.  Notice that, in the notation of Section~\ref{2.1},

\[ \begin{array}{rcl}

\widehat{P}(\rho)  &  =  &  \displaystyle \sum_{\tau \kappa \in T_2 K} P(\tau \kappa) \rho(\tau \kappa)
\ \  =  \ \ \displaystyle \frac{1}{r(n-r)r!(n-r)!} \sum_{\tau \in T_2} \sum_{\kappa \in K} \rho(\tau) \rho(\kappa) \vspace{1pc}  \\
&  =  &  \displaystyle \frac{1}{r(n-r)} \left[ \sum_{\tau \in T} \rho(\tau) \ - \ \sum_{\tau \in T_1} \rho(\tau) \right] \cdot
\left[ \frac{1}{r!(n-r)!} \sum_{\kappa \in K} \rho(\kappa) \right]. \vspace{.5pc}  \\

\end{array} \]

\noindent
Since $T$ is a conjugacy class in $S_n$, it follows from Lemma~\ref{2.4.1} that

\[ \sum_{\tau \in T} \rho(\tau) \ \ = \ \ |T| \ r(\rho) \ I \ \ = \ \ \mbox{$\frac{1}{2}$} n(n-1) \ r(\rho) \ I, \]

\noindent
where $r(\rho) := \chi_{\rho}(\tau) / d_{\rho}$, with $\tau \in T \subseteq S_n$, and $I$ is the $d_{\rho}$-dimensional identity
matrix. 

Since $L(X) \cong \rho_0 \uparrow_K^{S_n}$, where $\rho_0$ is the trivial representation of $K$, it follows that for any $\rho$ 
occurring in the decomposition of $L(X)$, $\rho \downarrow_K^{S_n}$ is the direct sum of $d_{\rho}$ copies of $\rho_0$.  Thus since 
$T_1 \subseteq K$, we have

\[ \sum_{\tau \in T_1} \rho(\tau) \ \ = \ \ |T_1| \ I \ \ = \ \ \left[ \mbox{$\frac{1}{2}$} n(n-1) \ - \ r(n-r) \right] I. \]

Now choose an orthonormal basis in $V$ such that the first $m_{\rho}$ basis vectors are $K$-invariant, as described in Section~\ref{2.2}.
It then follows from Lemma~\ref{2.2.3} that, in this basis, 

\[ \frac{1}{r!(n-r)!} \sum_{\kappa \in K} \rho(\kappa) \ \ = \ \ \widehat{I}(\rho). \]

\noindent
Therefore, it follows from these results that

\[ \widehat{P}(\rho) \ \ = \ \ \left[ \frac{\mbox{$\frac{1}{2}$} n(n-1) r(\rho) \ - \ \mbox{$\frac{1}{2}$} n(n-1) \ + \
r(n-r)}{r(n-r)} \right] \widehat{I}(\rho). \]

We now calculate $r(\rho)$ by use of Lemma~\ref{2.4.2}.  For the $i$th nontrivial irreducible representation $\rho = \rho_{[n-i,i]}$
of $S_n$ occurring in the decomposition of $L(X)$, it follows that, for $1 \leq i \leq r$,

\[ \begin{array}{rcl}

r(\rho)  &  =  &  \displaystyle \frac{1}{n(n-1)} \left[ (n-i)^2 \ - \ (n-i) \ + \ i^2 \ - \ 3i \right]  \vspace{1pc} \\
&  =  &  \displaystyle \frac{1}{n(n-1)} \left[ (n-i)(n-i-1) \ + \ i(i-3) \right]

\end{array} \]

\noindent
and hence

\[ \begin{array}{rcl}

\widehat{P}(\rho)  &  =  &  \displaystyle \left[ \frac{\mbox{$\frac{1}{2}$} (n-i)(n-i-1) \ + \ \mbox{$\frac{1}{2}$} i(i-3) \ - \
\mbox{$\frac{1}{2}$} n(n-1) \ + \ r(n-r)}{r(n-r)} \right] \widehat{I}(\rho)  \vspace{1pc} \\
&  =  &  \displaystyle \left[ 1 \ - \ \frac{i(n-i+1)}{r(n-r)} \right] \widehat{I}(\rho). \ \ \ \ \ \ \qed

\end{array} \]

We have now established the results necessary to prove Theorem~\ref{2.1.3}.  Recall from Section~\ref{2.2} that the probability 
measure $P$ defined in (\ref{2.1.1}) is bi-$K$-invariant.  By applying the result from Lemma~\ref{2.4.3} to the Upper Bound Lemma
(\ref{2.2.4}), we find that

\[ \begin{array}{rcl}

\displaystyle \| \widetilde{P^{*k}} - \widetilde{U} \|_{\mbox{\rm \scriptsize TV}}^2  &  \leq  &  \displaystyle \mbox{$\frac{1}{4}
{n \choose r}$} \| \widetilde{P^{*k}} - \widetilde{U} \|_2^2  \vspace{1pc} \\

&  =  &  \displaystyle \mbox{$\frac{1}{4}$} \sum_{\rho} d_{\rho} \ m_{\rho} \ \left[ 1 \ - \ \frac{i(n-i+1)}{r(n-r)} \right]^{2k},

\end{array} \]

\noindent
where the sum is taken over all nontrivial irreducible representations $\rho = \rho_{[n-i,i]}$ of $S_n$ occurring in the decomposition 
of $L(X)$.

Since, for each of these irreducible representations, $d_{\rho}$ was determined in Lemma~\ref{2.3.1} and $m_{\rho} = 1$, it then
follows that

\[ \begin{array}{rcl}

\displaystyle \| \widetilde{P^{*k}} - \widetilde{U} \|_{\mbox{\rm \scriptsize TV}}^2  &  \leq  &  \displaystyle \mbox{$\frac{1}{4} 
{n \choose r}$} \| \widetilde{P^{*k}} - \widetilde{U} \|_2^2  \vspace{1pc} \\

&  =  &  \displaystyle \mbox{$\frac{1}{4}$} \sum_{i=1}^r \left[ {n \choose i} \ - \ {n \choose i-1} \right] \cdot \left[ 1 \ - \ 
\frac{i(n-i+1)}{r(n-r)} \right]^{2k}  \vspace{1pc} \\

&  \leq  &  \displaystyle \mbox{$\frac{1}{4}$} \sum_{i=1}^r \left[ {n \choose i} \ - \ {n \choose i-1} \right] \exp \left\{ 
-\frac{2ki(n-i+1)}{r(n-r)} \right\}.

\end{array} \]

\noindent
Thus, when $k = \frac{1}{2} r \left( 1 - \frac{r}{n} \right) \left( \log n + c \right)$,

\[ \begin{array}{rcl}

\displaystyle \| \widetilde{P^{*k}} - \widetilde{U} \|_{\mbox{\rm \scriptsize TV}}^2  &  \leq  &  \displaystyle \mbox{$\frac{1}{4}
{n \choose r}$} \| \widetilde{P^{*k}} - \widetilde{U} \|_2^2  \vspace{1pc} \\

&  \leq  &  \displaystyle \mbox{$\frac{1}{4}$} \sum_{i=1}^r \left[ {n \choose i} \ - \ {n \choose i-1} \right] n^{-i(n-i+1)/n} 
e^{-ci(n-i+1)/n}.

\end{array} \]

A detailed discussion in Section 3 of Diaconis and Shahshahani (1987) determines the existence of a universal constant $a > 0$ such
that, for $c > 0$,
\num \begin{equation} \label{2.4.4}
\mbox{$\frac{1}{4}$} \sum_{i=1}^r \left[ {n \choose i} \ - \ {n \choose i-1} \right] n^{-i(n-i+1)/n} e^{-ci(n-i+1)/n} \ \ \leq \ \
a^2 e^{-c}.
\end{equation}

\noindent
This completes the proof of Theorem~\ref{2.1.3}. \vspace{.5pc}

Theorem~\ref{2.1.3} shows that $k = \frac{1}{2} r \left(1 - \frac{r}{n}\right) \left(\log n + c\right)$ steps are sufficient for the
(normalized) $\ell^2$ distance, and hence also the total variation distance, to become small.  Diaconis and Shahshahani (1987)
established a matching lower bound in the special case $r = n/2$.

\subsection{Analysis of a Variant of the Classical Model.} \label{2.5}

In preparation for our analysis in Section~\ref{3}, we now introduce a variant of the Bernoulli--Laplace diffusion model.
Suppose that the balls and racks are as described in Section~\ref{2.1}.  At each step of our variant of the classical model, two
positions $p$ and $q$ are chosen independently and uniformly from $\{1, 2, \ldots, n\}$.  If $p \neq q$, switch the balls in
positions $p$ and $q$.  If $p = q$ (which occurs with probability $1/n$), leave the balls in their current positions.  Then, if
necessary, permute the balls on each of the two racks so that their labels are in increasing order, as in the classical model.  

This switching procedure is exactly that of the random walk on the set of all possible positionings of $n$ cards introduced in Section
1 of Diaconis and Shahshahani (1981).  Furthermore, if $1 \leq p,q \leq r$ or $(r+1) \leq p,q \leq n$, then the resulting state of the
system is unchanged; this occurs with probability $1 - \frac{2r(n-r)}{n^2}$.  This procedure induces a Markov chain  on the same state
space as the classical model, but slows down the process by a factor of $\frac{n^2}{2r(n-r)}$ by not forcing two balls to be
switched at each step.  In light of Theorem~\ref{2.1.3}, one would certainly then expect $\frac{1}{4}n (\log n + c)$ steps to
suffice for near-stationarity.  We establish this conjecture as Theorem~\ref{2.5.3}.

Notice that $K \cong S_r~\times~S_{n-r}$ is also the stabilizer of the initial configuration of this process.  Let $T$, $T_1$ and
$T_2$ be defined as in Section \ref{2.1}.  As with the classical model, before the two positions $p$ and $q$ have been chosen, the
balls in the first rack may be permuted among themselves and the balls in the second rack may be permuted among themselves, without
changing the state of the system following the switch.  Thus, at each step, the procedure actually chooses a random element of $T_2
K := \{ \tau \kappa \in S_n  : \tau \in T_2$ \ and\ $\kappa \in K \}$ with probability $\frac{2r(n-r)}{n^2}$ and chooses
a random element of $K$ with probability $1 - \frac{2r(n-r)}{n^2}$.

As in Section~\ref{2.1}, this procedure may be modeled formally by a probability measure $Q$ on the symmetric group $S_n$.  Since
$Q(\tau\kappa) = \frac{2r(n-r)}{n^2} P(\tau\kappa)$ for all $\tau\kappa \in T_2 K$, where $P$ is defined in (\ref{2.1.1}), we thus
arrive at the following probability measure on the set of all permutations of $S_n$:
\num \begin{equation} \label{2.5.1}
\begin{array}{rcll}

Q(\kappa) & := & \displaystyle \frac{n^2 - 2r(n-r)}{n^2 r!(n-r)!} & \mbox{where $\kappa \in K$}, \vspace{.5pc} \\

Q(\tau\kappa) & := & \displaystyle \frac{2}{n^2r!(n-r)!} & \mbox{where $\tau \kappa \in T_2 K$, and}
\vspace{.5pc} \\

Q(\pi) & := & 0 & \mbox{otherwise}.

\end{array}
\end{equation}

For this probability measure $Q$, we now calculate the Fourier transform at each nontrivial irreducible representation of $S_n$
occurring in the decomposition of $L(X)$, as was done in Section~\ref{2.4} for $P$.

\begin{lemma} \label{2.5.2}
Let $Q$ be the probability measure on $S_n$ defined in (\ref{2.5.1}).  Let $\rho = \rho_{[n-i,i]}$ (for some $1 \leq i
\leq r$) be a nontrivial irreducible representation of $S_n$ occurring in the decomposition of $L(X)$.  Then, in a certain basis, 
the Fourier transform is

\[ \widehat{Q}(\rho) \ \ = \ \ \left[ 1 \ - \ \frac{2i(n-i+1)}{n^2} \right] \widehat{I}(\rho) \]

\noindent
where $\widehat{I}(\rho)$ is the $d_{\rho}~\times~d_{\rho}$ matrix whose only nonzero entry is a ``\ $1$'' in the $(1,1)$ position.
\end{lemma}

\proof{Proof} Let $P$ be the probability measure on $S_n$ defined in (\ref{2.1.1}).  Notice that

\[ Q \ \ = \ \ \left( 1 - \frac{2r(n-r)}{n^2} \right) U_K \ + \ \frac{2r(n-r)}{n^2} P, \]

\noindent
where $U_K$ is the probability measure on $S_n$ defined by $U_K(\kappa) := \frac{1}{|K|}$ for all $\kappa \in K$ and $U_K(\pi) := 0$
otherwise.  It then follows from Lemmas~\ref{2.2.3} and \ref{2.4.3}, using the same basis, that

\[ \begin{array}{rcl}

\widehat{Q}(\rho)  &  =  &  \displaystyle \sum_{\pi \in S_n} Q(\pi) \rho(\pi) \vspace{1pc} \\

&  =  &  \displaystyle \left( 1 - \frac{2r(n-r)}{n^2} \right) \sum_{\pi \in S_n} U_K(\pi) \rho(\pi) \ + \ \frac{2r(n-r)}{n^2}
\sum_{\pi \in S_n} P(\pi) \rho(\pi) \vspace{1pc} \\

&  =  &  \displaystyle \left( 1 - \frac{2r(n-r)}{n^2} \right) \widehat{I}(\rho) \ + \ \frac{2r(n-r)}{n^2} \widehat{P}(\rho) 
\vspace{1pc} \\

&  =  &  \displaystyle \left[ \left( 1 - \frac{2r(n-r)}{n^2} \right) \ + \ \frac{2r(n-r)}{n^2} \left( 1 \ - \
\frac{i(n-i+1)}{r(n-r)} \right) \right] \widehat{I}(\rho) \vspace{1pc} \\

&  =  &  \displaystyle \left[ 1 \ - \ \frac{2i(n-i+1)}{n^2} \right] \widehat{I}(\rho). \ \ \ \ \ \ \qed

\end{array} \]

The following result establishes an upper bound on both the total variation distance and the $\ell^2$ distance between 
$\widetilde{Q^{*k}}$ and $\widetilde{U}$, where $\widetilde{Q^{*k}}$ is the probability measure on the homogeneous space $X = S_n /
(S_r~\times~S_{n-r})$, induced by the convolution $Q^{*k}$ of $Q$ with itself $k$ times, and $\widetilde{U}$ is the uniform
probability measure on $X$.

\begin{theorem} \label{2.5.3}
Let $Q$ and $U$ be the probability measures on the symmetric group $S_n$ defined in (\ref{2.5.1}) and (\ref{2.1.2}),
respectively.  Let $\widetilde{Q^{*k}}$ be the probability measure on the homogeneous space $X = S_n / (S_r~\times~S_{n-r})$ induced
by $Q^{*k}$ and let $\widetilde{U}$ be the uniform probability measure on $X$.  Let $k = \frac{1}{4} n (\log n + c)$.  Then there
exists a universal constant $a > 0$ such that

\[ \| \widetilde{Q^{*k}} - \widetilde{U} \|_{\mbox{\rm \scriptsize TV}} \ \ \leq \ \ \mbox{$\frac{1}{2} {n \choose r}^{1/2}$} \|
\widetilde{Q^{*k}} - \widetilde{U} \|_2 \ \ \leq \ \ ae^{-c/2} \ \ \ \mathrm{for\ all\ } \mbox{$c > 0$}. \]

\end{theorem}

\noindent
This Markov chain is twice as fast to converge as the random walk on the symmetric group $S_n$ introduced in Section 1 of Diaconis and
Shahshahani (1981).

\proof{Proof} Notice that $Q$ is clearly right $K$-invariant.  (In fact, it is bi-$K$-invariant.)  By applying the result of
Lemma~\ref{2.5.2} to the Upper Bound Lemma (\ref{2.2.4}), we find that

\[ \displaystyle \| \widetilde{Q^{*k}} - \widetilde{U} \|_{\mbox{\rm \scriptsize TV}}^2  \ \ \leq \ \ \displaystyle \mbox{$\frac{1}{4}
{n \choose r}$} \| \widetilde{Q^{*k}} - \widetilde{U} \|_2^2
\ \ = \ \ \displaystyle \mbox{$\frac{1}{4}$} \sum_{\rho} d_{\rho} \ m_{\rho} \ \left[ 1 \ - \ \frac{2i(n-i+1)}{n^2} \right]^{2k}, \]

\noindent
where the sum is taken over all nontrivial irreducible representations $\rho = \rho_{[n-i,i]}$ of $S_n$ occurring in the decomposition 
of $L(X)$.

Since, for each of these irreducible representations, $d_{\rho}$ was determined in Lemma~\ref{2.3.1} and $m_{\rho} = 1$, it then
follows that
\num \begin{equation} \label{2.5.4}
\begin{array}{rcl}

\displaystyle \| \widetilde{Q^{*k}} - \widetilde{U} \|_{\mbox{\rm \scriptsize TV}}^2  &  \leq  &  \displaystyle \mbox{$\frac{1}{4}
{n \choose r}$} \| \widetilde{Q^{*k}} - \widetilde{U} \|_2^2  \vspace{1pc} \\

&  =  &  \displaystyle \mbox{$\frac{1}{4}$} \sum_{i=1}^r \left[ {n \choose i} \ - \ {n \choose i-1} \right] \cdot \left[ 1 \ - \ 
\frac{2i(n-i+1)}{n^2} \right]^{2k}  \vspace{1pc} \\

&  \leq  &  \displaystyle \mbox{$\frac{1}{4}$} \sum_{i=1}^r \left[ {n \choose i} \ - \ {n \choose i-1} \right] \exp \left\{ 
-\frac{4ki(n-i+1)}{n^2} \right\}.  \vspace{1pc} \\

\end{array}
\end{equation}

\noindent
Thus, when $k = \frac{1}{4} n (\log n + c)$,

\[ \begin{array}{rcl}

\displaystyle \| \widetilde{Q^{*k}} - \widetilde{U} \|_{\mbox{\rm \scriptsize TV}}^2  &  \leq  &  \displaystyle \mbox{$\frac{1}{4}
{n \choose r}$} \| \widetilde{Q^{*k}} - \widetilde{U} \|_2^2  \vspace{1pc} \\

&  \leq  &  \displaystyle \mbox{$\frac{1}{4}$} \sum_{i=1}^r \left[ {n \choose i} \ - \ {n \choose i-1} \right] n^{-i(n-i+1)/n} 
e^{-ci(n-i+1)/n}.

\end{array} \]

\noindent
It then follows from (\ref{2.4.4}) that, for $c > 0$,

\[ \mbox{$\frac{1}{4}$} \sum_{i=1}^r \left[ {n \choose i} \ - \ {n \choose i-1} \right] n^{-i(n-i+1)/n} e^{-ci(n-i+1)/n} \ \ \leq \
\ a^2 e^{-c}, \]

\noindent
from which the desired result follows. \qed

Theorem~\ref{2.5.3} shows that $k = \frac{1}{4} n \left(\log n + c\right)$ steps are sufficient for the (normalized) $\ell^2$ 
distance, and hence the total variation distance, to become small.  A lower bound in the (normalized) $\ell^2$ metric can also be
derived by examining   $(n-1) \left( 1 - \frac{2}{n} \right)^{2k}$, which is the dominant contribution to the summation
(\ref{2.5.4}) from the proof of Theorem~\ref{2.5.3}.  This term corresponds to the choice $i = 1$.  Notice that $k = \frac{1}{4} n
\left( \log n - c \right)$ steps are necessary for just this term to become small.

A number of steps necessary for the total variation distance to become small is given by the following result.  An analogous result
was proved for the classical model (\ref{2.1.1}), in the special case $r = n/2$, in Theorem 1 of Diaconis and Shahshahani (1987).

\begin{theorem} \label{2.5.5}
Let $Q$ and $U$ be the probability measures on the symmetric group $S_n$ defined in (\ref{2.5.1}) and (\ref{2.1.2}), respectively.
Let $\widetilde{Q^{*k}}$ be the probability measure on the homogeneous space $X = S_n / (S_r~\times~S_{n-r})$ induced
by $Q^{*k}$ and let $\widetilde{U}$ be the uniform probability measure on $X$.  Let $n \geq 3$ and 

\[ \mbox{$k = \min \left\{ \frac{1}{4} n (\log n - c), \ (n-2) \log \left( \frac{n}{n-2r} \right) + \frac{1}{2} (n-2) \log \left[ 1
+ \frac{1}{4} (n-2) \hspace{-.05in} \left( 1 - \left( \frac{n-2r}{n} \right)^2 \right) e^{-c} \right] \right\}$} \]

\noindent
be a nonnegative integer, with $c \geq 0$ (and where $\log \left( \frac{n}{0} \right)$ is interpreted as $+\infty$ when $r = n/2$).  
Then there exists a universal constant $\tilde{a} > 0$ such that

\[ \mbox{$\frac{1}{2} {n \choose r}^{1/2}$} \| \widetilde{Q^{*k}} - \widetilde{U} \|_2 \ \ \geq \ \ \|
\widetilde{Q^{*k}} - \widetilde{U} \|_{\mbox{\rm \scriptsize TV}} \ \ \geq \ \ 1 - \tilde{a}e^{-c}. \]

\end{theorem}

\proof{Proof} Recall from Lemma~\ref{2.3.1} that for $L(X)$, the set of all (complex-valued) functionals defined on $X = S_n / (S_r
\times S_{n-r})$, we have the multiplicity-free decomposition

\[ \displaystyle L(X) \ \ \cong \ \ \rho_{[n]} \ \oplus \ \rho_{[n-1,1]} \ \oplus \ \cdots \ \oplus \ \rho_{[n-r,r]}, \]

\noindent
where $\rho_{[n-i,i]} : S_n \longrightarrow$ \textbf{GL} $(V_i)$ (say).  So it follows from Lemma~\ref{2.2.2} that each $V_i$, for
$0 \leq i \leq r$, has a unique nontrivial one-dimensional space of $K$-fixed vectors, where $K := S_r~\times~S_{n-r}$.  

The \emph{spherical function} $s_i$, for $0 \leq i \leq r$, is defined to be the unique left $K$-invariant function (i.e., vector) 
in $V_i$ normalized so that $s_i(x_0) = 1$, where $x_0$ is the left coset representative corresponding to $e \in S_n$.  It follows
from Exercise 17 in Section F of Chapter 3 of Diaconis (1988) that the spherical functions satisfy

\[ s_i(x) \ = \ \frac{1}{|K|} \sum_{\kappa \in K} \chi_i(x \kappa) \ \ \ \mathrm{for\ } \mbox{$x \in X$}, \]

\noindent
where $\chi_i$ is the character of $\rho_{[n-i,i]}$, with $[n-0,0] \equiv [n]$.

Under the uniform measure $\widetilde{U}$ on $X$, notice that

\[ E_{\widetilde{U}}( s_i ) \ \ = \ \ \frac{1}{|X|} \sum_{x \in X} s_i(x) \ \ = \ \ \frac{1}{|X|}
\sum_{x \in X} \frac{1}{|K|} \sum_{\kappa \in K} \chi_i(x \kappa) \ \ = \ \ \frac{1}{n!} \sum_{\pi \in S_n} \chi_i(\pi) \ \ = \ \
\langle \chi_i, \chi_0 \rangle_{S_n}, \]

\noindent
where $\chi_0$ is the character of the trivial representation $\rho_{[n]}$.  Thus since the irreducible characters of a group $G$ form
an orthonormal basis for the Hilbert space of class functions on $G$, it then follows that $E_{\widetilde{U}}( s_0 ) = 1$ and
$E_{\widetilde{U}}( s_i ) = 0$ for $1 \leq i \leq r$.

Under the $k$-fold convolution measure $\widetilde{Q^{*k}}$ on $X$, notice that

\[ \begin{array}{rcl}

\displaystyle E_{\widetilde{Q^{*k}}}( s_i )  &  =  &  \displaystyle \sum_{x \in X} \widetilde{Q^{*k}}(x) s_i(x) \ = \ \sum_{x \in X}
\widetilde{Q^{*k}}(x) \frac{1}{|K|} \sum_{\kappa \in K} \chi_i(x \kappa) \vspace{1pc} \\

&  =  &  \displaystyle \sum_{\pi \in S_n} Q^{*k}(\pi) \chi_i(\pi) \ \ = \ \ \mathrm{tr} \mbox{$\displaystyle \sum_{\pi \in S_n}
Q^{*k}(\pi) \rho_{[n-i,i]} (\pi) \ \ = \ \ $} \mathrm{tr} \ \mbox{$\widehat{Q^{*k}}(\rho_{[n-i,i]})$},

\end{array} \]

\noindent
where we use the fact that $\widetilde{Q^{*k}}(x) = |K| \ Q^{*k} (x \kappa)$ for each $\kappa \in K$.  In particular, it follows
from Lemma~\ref{2.5.2} that

\[ E_{\widetilde{Q^{*k}}}( s_1 ) \ = \ \left( 1 - \mbox{$\frac{2}{n}$} \right)^k \ \ \ \mathrm{and} \ \ \
\mbox{$E_{\widetilde{Q^{*k}}}( s_2 ) \ = \ \left( 1 - \mbox{$\frac{2}{n}$} \right)^{2k}$}. \]

Define $f(x) := \sqrt{n-1} \ s_1(x)$ for $x \in X$.  Then

\[ E_{\widetilde{Q^{*k}}}( f ) \ = \ \sqrt{n-1} \left( 1 - \mbox{$\frac{2}{n}$} \right)^k \ \ \ \mathrm{and} \ \ \ 
\mbox{$E_{\widetilde{U}}( f ) \ = \ 0$}. \]

\noindent
In order to determine Var$_{\widetilde{U}}( f )$ and Var$_{\widetilde{Q^{*k}}}( f )$, we must calculate $E_{\widetilde{U}}( f^2 )$
and $E_{\widetilde{Q^{*k}}}( f^2 )$.  This is done with the aid of the following identity (which can be derived from formulas for
$s_1$ and $s_2$, as suggested in Diaconis and Shahshahani (1987)):

\[ s_1^2 \ \ = \ \ \frac{1}{n-1} \ + \frac{(n-2r)^2}{r(n-r)(n-2)} s_1 \ + \ \frac{n^2(r-1)(n-r-1)}{r(n-r)(n-1)(n-2)} s_2. \]

\noindent
Thus

\[ \begin{array}{rcl}

\mbox{Var}_{\widetilde{U}}( f )  &  =  &  \displaystyle (n-1) E_{\widetilde{U}}( s_1^2 ) \vspace{1pc} \\ 

&  =  &  \frac{n-1}{n-1} \ + \ \frac{(n-1)(n-2r)^2}{r(n-r)(n-2)} E_{\widetilde{U}}( s_1 ) \ + \ 
\frac{(n-1)n^2(r-1)(n-r-1)}{r(n-r)(n-1)(n-2)} E_{\widetilde{U}}( s_2 ) \ \ = \ \ 1

\end{array} \]

\noindent
and

\[ \begin{array}{rcl}

\mbox{Var}_{\widetilde{Q^{*k}}}( f )  &  =  &  \displaystyle (n-1) E_{\widetilde{Q^{*k}}}( s_1^2 ) \ - \ (n-1) \left(
E_{\widetilde{Q^{*k}}}( s_1 ) \right)^2 \vspace{1pc} \\

&  =  &  \frac{n-1}{n-1} \ + \ \frac{(n-1)(n-2r)^2}{r(n-r)(n-2)} \left( 1 - \mbox{$\frac{2}{n}$} \right)^k \vspace{1pc} \\

&     &  + \ \ \frac{(n-1)n^2(r-1)(n-r-1)}{r(n-r)(n-1)(n-2)} \left( 1 - \mbox{$\frac{2}{n}$} \right)^{2k} \ - \ (n-1) \left( 1 -  
\mbox{$\frac{2}{n}$} \right)^{2k} \vspace{1pc} \\

&  =  &  1 \ + \ \frac{4\left(\frac{n-1}{n-2}\right) (n-2r)^2} {n^2 - (n-2r)^2} \left( 1 - \mbox{$\frac{2}{n}$}
\right)^k \ - \ \left[ \frac{4\left(\frac{n-1}{n-2}\right) n^2} {n^2 - (n-2r)^2} - \frac{3n-2}{n-2} \right] \left( 1 -
\mbox{$\frac{2}{n}$} \right)^{2k}.

\end{array} \] 

\noindent
Therefore,
\num \begin{eqnarray} \label{2.5.6} 
\displaystyle \frac{\mbox{Var}_{\widetilde{Q^{*k}}}(f)}{\left(E_{\widetilde{Q^{*k}}}(f)\right)^2}  &  =  &  \displaystyle
\frac{1}{\left(E_{\widetilde{Q^{*k}}}(f)\right)^2} \vspace{1pc} \\
&  +  &  \displaystyle \frac{\frac{4n^2}{n-2}}{n^2 - (n-2r)^2} \left[ \left( \frac{n-2r}{n} \right)^2 \left( 1 - \frac{2}{n}
\right)^{-k} - 1 \right] \ \ + \ \ \frac{3n-2}{(n-1)(n-2)}. \nonumber
\end{eqnarray}

By elementary calculus, $x \leq -\log(1-x) \leq \frac{x}{1-x}$ for $0 \leq x \leq 1$.  Thus if $k \leq \frac{1}{4} n \left( \log n -
c \right)$, with $n \geq 3$ and $c \geq 0$, then

\[ \begin{array}{rcl}

E_{\widetilde{Q^{*k}}}( f )  &  =  &  \sqrt{n-1} \left( 1 - \mbox{$\frac{2}{n}$} \right)^k \ \ \geq \ \ \sqrt{n-1} \ e^{-2k/(n-2)}
\vspace{1pc} \\

&  \geq  &  \left( 1 - \mbox{$\frac{1}{n}$} \right)^{1/2} \left( \frac{1}{n} \right)^{1/(n-2)} e^{c/2} e^{c/(n-2)} \ \ \geq \ \
\sqrt{\frac{2}{27}} \ e^{c/2},

\end{array} \]

\noindent
where we note that, for $n \geq 3$, $\left( 1 - \mbox{$\frac{1}{n}$} \right)^{1/2} \left( \frac{1}{n} \right)^{1/(n-2)}$ is
increasing and $e^{c/(n-2)} \geq 1$.

Also notice that the first and third terms on the right in (\ref{2.5.6}) are bounded by

\[ \begin{array}{rcl}

\displaystyle \frac{1}{\left(E_{\widetilde{Q^{*k}}}(f)\right)^2}  &  \leq  &  \frac{27}{2} e^{-c} \ \ \ \mathrm{and} \vspace{1pc} \\

\displaystyle \frac{3n-2}{(n-1)(n-2)}  &  \leq  &  \displaystyle \frac{3n}{\left(\frac{2}{3}n\right) \left(\frac{1}{3}n\right)} \ \
\leq \ \ \mbox{$\frac{27}{2}$} n^{-1} \ \ \leq \ \ \mbox{$\frac{27}{2}$} e^{-c}

\end{array} \]

\noindent
when $n \geq 3$ and $0 \leq c \leq \log n$.  In order for the second term on the right in (\ref{2.5.6}) to be bounded above by
$e^{-c}$ we must have
\num \begin{equation} \label{2.5.7}
\left( \frac{n-2r}{n} \right)^2 \left( 1 - \frac{2}{n} \right)^{-k}\ \ \leq \ \ 1 \ + \ \frac{(n-2)[n^2 - (n-2r)^2]}{4n^2} \ e^{-c}.
\end{equation}

\noindent
Since, when $n \geq 3$,

\[ \left( 1 - \frac{2}{n} \right)^{-k} \ \ \leq \ \ e^{2k/(n-2)}, \]

\noindent
it is sufficient for (\ref{2.5.7}) to have

\[ \frac{2k}{n-2} \ \ \leq \ \ 2 \log \left( \frac{n}{n-2r} \right) \ + \ \log \left[ 1 \ + \ \frac{(n-2)[n^2 - (n-2r)^2]}{4n^2} \
e^{-c} \right], \]

\noindent
i.e., to have

\[ k \ \ \leq \ \ (n-2) \log \left( \frac{n}{n-2r} \right) \ + \ \mbox{$\frac{1}{2}$} (n-2) \log \left[ 1 \ + \ \frac{(n-2)[n^2 -
(n-2r)^2]}{4n^2} \ e^{-c} \right]. \]

\noindent
In summary, if $n \geq 3$ and $0 < c \leq \log n$, and if

\[ \mbox{$k = \min \left\{ \frac{1}{4} n (\log n - c), (n-2) \log \left( \frac{n}{n-2r} \right) + \frac{1}{2} (n-2) \log \left[ 1
+ \frac{1}{4} (n-2) \hspace{-.05in} \left( 1 - \left( \frac{n-2r}{n} \right)^2 \right) e^{-c} \right] \right\}$}, \]

\noindent
then

\[ E_{\widetilde{Q^{*k}}}( f ) \ \geq \ \sqrt{\mbox{$\frac{2}{27}$}} \ e^{c/2} \ \ \ \mathrm{and} \ \ \ 
\mbox{$\displaystyle \frac{\mbox{Var}_{\widetilde{Q^{*k}}}(f)} {\left(E_{\widetilde{Q^{*k}}}(f)\right)^2} \ \leq \ 28
e^{-c}$}. \]

Now define $A_{\alpha} := \{ x \in X : |f(x)| \leq \alpha \}$.  It follows from Chebyshev's inequality that $\displaystyle
\widetilde{U}(A_{\alpha}) \geq 1 - \frac{1}{\alpha^2}$ and that $\displaystyle \widetilde{Q^{*k}}(A_{\alpha}) \leq \frac{28e^{-c}} 
{\left( \mbox{$\sqrt{\frac{2}{27}}$} e^{c/2} - \alpha \right)^2}$, provided $0 \leq \alpha <  \mbox{$\sqrt{\frac{2}{27}}$} e^{c/2}$.
Then

\[ \mbox{$\frac{1}{2} {n \choose r}^{1/2}$} \| \widetilde{Q^{*k}} - \widetilde{U} \|_2 \ \ \geq \ \ \|
\widetilde{Q^{*k}} - \widetilde{U} \|_{\mbox{\rm \scriptsize TV}} \ \ \geq \ \ 1 \ - \ \frac{1}{\alpha^2} \ - \ \frac{28e^{-c}}{\left( 
\mbox{$\sqrt{\frac{2}{27}}$} \ e^{c/2} - \alpha \right)^2}. \]

\noindent
Choosing $\alpha = \mbox{$\frac{1}{2} \sqrt{\frac{2}{27}}$} \ e^{c/2}$ shows that

\[ \mbox{$\frac{1}{2} {n \choose r}^{1/2}$} \| \widetilde{Q^{*k}} - \widetilde{U} \|_2 \ \ \geq \ \ \|
\widetilde{Q^{*k}} - \widetilde{U} \|_{\mbox{\rm \scriptsize TV}} \ \ \geq \ \ 1 \ - \ 54 e^{-c} \ - \ 1512 e^{-2c} \ \ \geq \ \ 1 \ 
- \ 1566 e^{-c}, \]

\noindent
which completes the proof. \qed

Theorem~\ref{2.5.5} gives a number of steps necessary for the total variation distance to become small in our variant of the
classical Bernoulli--Laplace diffusion model.  We now examine in rough terms how the value of $r$ dictates the choice of $k$.  For
notational purposes, let

\[ \mbox{$ f(n,r,c) \ \ := \ \ (n-2) \log \left( \frac{n}{n-2r} \right) + \frac{1}{2} (n-2) \log \left[ 1 + \frac{1}{4} (n-2) \left(
1 - \left( \frac{n-2r}{n} \right)^2 \right) e^{-c} \right]$}. \]

\noindent
Notice that $f(n,r,c)$ is increasing in $r$.  At one extreme, let $r = n/2$ (assuming for simplicity that $n$ is even); then, since
$f(n,r,c) = +\infty$, $k$ is chosen to be $\frac{1}{4} n \left( \log n - c \right)$, matching the upper bound in
Theorem~\ref{2.5.3}.  At the other extreme, let $r = 1$; then it can be shown that $\frac{1}{33} e^{-c} n \leq f(n,r,c) \leq
\frac{1}{4} n \left( \log n - c \right)$, and so $k$ is chosen to be $f(n,r,c) \geq \frac{1}{33} e^{-c} n$.  In this case we find
that order $n$ (\emph{not} order $n \log n$) steps are necessary, and indeed it is easy to show that order $n$ steps are also 
sufficient.  There is some value of $r$ for which $f(n,r,c)$ ``crosses over'' $\frac{1}{4} n \left( \log n - c\right)$.  This occurs 
in the vicinity of $2r = n^{1/2}$.  But even for $2r$ as small as $n^{\delta}$, for any fixed $\delta > 0$, the value of $k$ 
determined in Theorem~\ref{2.5.5} is of order $n \log n$.

\section{A Signed Generalization of the Bernoulli--Laplace Diffusion Model.} \label{3}

\subsection{Introduction.} \label{3.1}

We now extend (the variant in Section~\ref{2.5} of) the Bernoulli--Laplace diffusion model to the case in which the balls
also have charges (positive or negative).  Imagine two racks, the first with positions labeled $1$ through $r$ and the second with
positions labeled $r+1$ through $n$.  Without loss of generality, we assume that $1 \leq r \leq n/2$.  Suppose that there are $n$
balls, labeled $1$ through $n$, each initially placed at its corresponding position.  Also suppose that each ball has a charge
(positive or negative) and that initially each ball is positively charged.  We refer to this as the \emph{signed} Bernoulli--Laplace
diffusion model.  

At each step, independently choose two positions $p$ and $q$ uniformly from $\{1, 2, \ldots,   n\}$.  If $p \neq q$, switch
the balls in positions $p$ and $q$.  Then independently, with probability $1/2$, change the charge of the ball moved to position
$p$; and independently, also with probability $1/2$, change the charge of the ball moved to position $q$.  Then, if necessary,
permute the balls on each of the two racks so that their labels are in increasing order.  If $p = q$ (which occurs with probability
$1/n$), leave the balls in their current positions.  Then, again independently with probability $1/2$, change the charge of the ball
in position $p = q$.

We refer to the process described above as the \emph{independent flips} model.  A similar process, known as the \emph{paired flips}
model, is introduced in Section~\ref{3.4}.

This switching procedure is exactly that of the random walk that was introduced in Section 3.1 of Schoolfield (1999) in the special
case of the hyperoctahedral group $\mathbb{Z}_2~\wr~S_n$.  The state of our signed Bernoulli--Laplace system is completely determined,
at each step, by the ordered $n$-tuple of charges of the $n$ balls $1, 2, \ldots, n$ and the labels of balls on the first rack.  Our
switching procedure thus induces a Markov chain on the set of all $2^n \cdot {n \choose r}$ ordered pairs of $n$-dimensional binary
vectors and $r$-element subsets of a set with $n$ elements.  

Let $K$ be the subgroup of $\mathbb{Z}_2~\wr~S_n$ which permutes the first $r$ indices among themselves and the last $n-r$ indices
among themselves, but does not commingle these two sets of indices.  Thus $K \cong S_r~\times~S_{n-r}$.  Notice that $K$ is the
stabilizer of the initial configuration of the process.  The switching procedure described above thus also induces a Markov chain on
the homogeneous space $X = (\mathbb{Z}_2~\wr~S_n) / (S_r~\times~S_{n-r})$.  The set $X$ may thus be identified with the set
of all $2^n \cdot {n \choose r}$ ordered pairs of $n$-dimensional binary vectors and $r$-element subsets of a set with $n$ elements.  

Let $T$ be the set of all signed transpositions in $\mathbb{Z}_2~\wr~S_n$.  Also let $T_1$ be the set of all signed transpositions 
in $K$, let $T_2$ be the set of all signed transpositions in $(\mathbb{Z}_2~\wr~K) \setminus K$, and let $T_3$ be the set of
all signed transpositions in $(\mathbb{Z}_2~\wr~S_n) \setminus (\mathbb{Z}_2~\wr~K)$.  Thus $T = T_1 \cup T_2 \cup
T_3$.  Notice that $\vec{v} = \vec{0} \in \mathbb{Z}_2^n$ for any $(\vec{v};\tau) \in T_1$, that $\vec{v} \in
\mathbb{Z}_2^n$ has one or two nonzero entries for any $(\vec{v};\tau) \in T_2$, and that $\vec{v} \in \mathbb{Z}_2^n$
has zero, one, or two nonzero entries for any $(\vec{v};\tau) \in T_3$.  Finally, let $U$ be the set of all signed identities in
$\mathbb{Z}_2~\wr~S_n$.  Notice that $U \subseteq \mathbb{Z}_2~\wr~K$.  Recall that for any $(\vec{u};e) \in U$, $\vec{u} \in
\mathbb{Z}_2^n$ has exactly one nonzero entry.

As with the classical model, before the two positions $p$ and $q$ have been chosen, the balls in the first rack may be permuted  
among themselves and the balls in the second rack may be permuted among themselves, without changing the state of the system
following the switch.  Thus, at each step, a random element of $K$ is effectively generated whenever the procedure described above
results in the identity or in a signed transposition in $T_1$; this event occurs with probability 

\[ \frac{1}{2n} + \frac{\frac{1}{2} n(n-1) - r(n-r)}{2n^2} = \frac{n(n+1) - 2r(n-r)}{4n^2}. \]

A similar analysis shows that the procedure effectively generates a random element of $U K = \{ (\vec{u};\kappa) \in
\mathbb{Z}_2~\wr~S_n : (\vec{u};e) \in U \mbox{\ and\ } (\vec{0};\kappa) \in K \}$ with probability $\frac{1}{2n}$, a
random element of $T_2 K = \{ (\vec{v};\tau\kappa) \in \mathbb{Z}_2~\wr~S_n : (\vec{v};\tau) \in T_2 \mbox{\ and\ }
(\vec{0};\kappa) \in K \}$ with probability $\frac{3n(n-1) - 6r(n-r)}{4n^2}$, and a random element of $T_3 K = \{
(\vec{v};\tau\kappa) \in \mathbb{Z}_2~\wr~S_n : (\vec{v};\tau) \in T_3 \mbox{\ and\ } (\vec{0};\kappa) \in K \}$ with
probability $\frac{2r(n-r)}{n^2}$.

Notice that each element of $T_3 K$ can be uniquely written as $(\vec{v};\tau\kappa)$, where $(\vec{v};\tau) \in T_3$ and 
$(\vec{0};\kappa) \in K$.  However, $U K \subseteq T_2 K$, with the exception that the elements $(\vec{u};\kappa) \in U K$ with 
$\vec{u} = (1, 0, \ldots, 0)$ are not included in $T_2 K$ when $r = 1$.  But each element of $U K$ can be uniquely written as
$(\vec{u};\kappa)$, where $(\vec{u};e) \in U$ and $(\vec{0};\kappa) \in K$ and each element of $T_2 K \setminus U K$ (where the set
difference here is proper unless $r = 1$) can be uniquely written as $(\vec{v};\tau\tau^{-1}\kappa)$, where $(\vec{v};\tau) \in T_2$
and $(\vec{0};\tau^{-1}\kappa) \in K$.  Let $U_r$ consist of the signed identities of the first $r$ indices and $U_{n-r}$ consist of
the signed identities of the last $n-r$ indices.  Notice that 

\[ \begin{array}{rcl}

|K|  &  =  &  r!(n-r)!, \vspace{.5pc} \\
|U_r K|  &  =  &  r \, r!(n-r)!, \vspace{.5pc} \\
|U_{n-r} K|  &  =  &  (n-r)r!(n-r)!, \vspace{.5pc} \\
|T_2 K \setminus U K|  &  =  &  \left[ \mbox{$\frac{1}{2}$} n(n-1) - r(n-r) \right] r!(n-r)!, \mathrm{\ and} \vspace{.5pc} \\
|T_3 K|  &  =  &  4r(n-r) \cdot r!(n-r)!.

\end{array} \]

The signed Bernoulli--Laplace diffusion model may be modeled formally by a probability measure $P$ on the hyperoctahedral group
$\mathbb{Z}_2~\wr~S_n$.  We may thus define the following probability measure on the set of all signed permutations of $\mathbb{Z}_2
\wr S_n$:
\num \begin{equation} \label{3.1.1}
\begin{array}{rcll}

P(\vec{0};\kappa) & := & \frac{n(n+1) - 2r(n-r)}{4n^2 r!(n-r)!} & \mbox{where $(\vec{0};\kappa)
\in K$}, \vspace{.5pc} \\

P(\vec{u};\kappa) & := & \frac{r}{2n^2 r!(n-r)!} & \mbox{where $(\vec{u};\kappa) \in U_r K$}, \vspace{.5pc} \\

P(\vec{u};\kappa) & := & \frac{n-r}{2n^2 r!(n-r)!} & \mbox{where $(\vec{u};\kappa) \in U_{n-r} K$}, \vspace{.5pc} \\

P(\vec{v};\tau\tau^{-1}\kappa) & := & \frac{1}{2n^2 r!(n-r)!} & \mbox{where $(\vec{v};\tau\tau^{-1}\kappa) \in T_2 K
\setminus U K$}, \vspace{.5pc} \\

P(\vec{v};\tau\kappa) & := & \frac{1}{2n^2 r!(n-r)!} & \mbox{where $(\vec{v};\tau\kappa) \in T_3 K$, and}
\vspace{.5pc} \\

P(\vec{x};\pi) & := & 0 & \mbox{otherwise}.

\end{array}
\end{equation}

\noindent
Since there are $2^n \cdot n!$ elements in $\mathbb{Z}_2~\wr~S_n$, the uniform probability measure on the set of all signed
permutations is given by
\num \begin{equation} \label{3.1.2}
U(\vec{x};\pi) \ \ := \ \ \mbox{$\frac{1}{2^n \cdot n!}$} \ \ \ \mbox{for every $(\vec{x};\pi) \in \mathbb{Z}_2~\wr~S_n$}. 
\end{equation}

The following result establishes an upper bound on both the total variation distance and the $\ell^2$ distance between
$\widetilde{P^{*k}}$ and $\widetilde{U}$, where $\widetilde{P^{*k}}$ is the probability measure on the homogeneous space $X =
(\mathbb{Z}_2~\wr~S_n) / (S_r~\times~S_{n-r})$ induced by the convolution $P^{*k}$ of $P$ with itself $k$ times, and
$\widetilde{U}$ is the uniform probability measure on $X$.  It should be noted that $\| \widetilde{P^{*k}} - \widetilde{U}
\|_{\mbox{\rm \scriptsize TV}}$ is the total variation distance to uniformity after $k$ steps for the Markov chain on $X$ defined 
above, started at the chosen initial configuration.  We establish an analogous result for the paired flips model as Theorem~\ref{3.4.3}.

\begin{theorem} \label{3.1.3}
Let $P$ and $U$ be the probability measures on the hyperoctahedral group $\mathbb{Z}_2~\wr~S_n$ defined in
(\ref{3.1.1}) and (\ref{3.1.2}), respectively.  Let $\widetilde{P^{*k}}$ be the probability measure on the homogeneous space $X =
(\mathbb{Z}_2~\wr~S_n) / (S_r~\times~S_{n-r})$ induced by $P^{*k}$ and let $\widetilde{U}$ be the uniform probability
measure defined on $X$.  Let $k = \frac{1}{4} n (\log n + c)$.  Then there exists a universal constant $b > 0$ such that

\[ \| \widetilde{P^{*k}} - \widetilde{U} \|_{\mbox{\rm \scriptsize TV}} \ \ \leq \ \ \mbox{$\frac{1}{2} \left[ 2^n \cdot {n \choose
r} \right]^{1/2}$} \| \widetilde{P^{*k}} - \widetilde{U} \|_2 \ \ \leq \ \ be^{-c/2} \ \ \ \mathrm{for\ all\ } 
\mbox{$c > 0$}.
\]

\end{theorem}

\noindent
Notice that this is (essentially) the same result as that found in Theorem~\ref{2.5.3}.  This Markov chain is twice as fast to
converge as the random walk analyzed in Section 3 of Schoolfield (1999) in the special case of the hyperoctahedral group $\mathbb{Z}_2
\wr S_n$.

In the following sections we present the results needed to prove this theorem and an analogous theorem for the paired flips model.
In Section~\ref{3.2} we decompose the representation $L(X)$ of $\mathbb{Z}_2~\wr~S_n$ into its irreducible constituents.  In
Section~\ref{3.3} we calculate the Fourier transform of the probability measure defined in (\ref{3.1.1}), and this is followed by
the proof of Theorem~\ref{3.1.3}.  In Section~\ref{3.4} we perform a similar analysis of the paired flips model.

\subsection{Irreducible Representations in $L(X)$.} \label{3.2}

The decomposition of the representation $L(X)$ of $\mathbb{Z}_2~\wr~S_n$ into its irreducible constituents is given by the 
following.  See Section 3.4 of Schoolfield (1999) for details about the irreducible representations of $G~\wr~S_n$ for any $G$.

\begin{lemma} \label{3.2.1}
Let $L(X)$ be the set of all (complex-valued) functionals defined on $X = (\mathbb{Z}_2~\wr~S_n) /
(S_r~\times~S_{n-r})$.  Then, as a representation of $\mathbb{Z}_2~\wr~S_n$,

\[ L(X) \ \ \cong \ \ \bigoplus_{i=0}^r \ \bigoplus_{j=0}^{n-r} \ \bigoplus_{\ell=0}^{i \wedge j} \ \bigoplus_{m=0}^{(r-i) \wedge
((n-r)-j)} \rho_{([(i+j)-\ell,\ell]; [n-(i+j)-m,m])} \]

\noindent
where $\rho_{([(i+j)-\ell,\ell]; [n-(i+j)-m,m])}$ is the irreducible representation of $\mathbb{Z}_2~\wr~S_n$
corresponding to the two-part partitions $[(i+j)-\ell,\ell]$ and $[n-(i+j)-m,m]$ of $i+j$ and $n-(i+j)$, respectively.
\end{lemma}

\proof{Proof} Recall that as a representation of $\mathbb{Z}_2~\wr~S_n$, $L(X)$ is isomorphic to $\rho_0 \uparrow_{S_r \times
S_{n-r}}^{\mathbb{Z}_2~\wr~S_n}$, where $\rho_0$ is the trivial representation of $S_r~\times~S_{n-r}$.  It follows from Theorem 10 
in Section 3.2 of Serre (1977) that $\rho_0 = \rho_0^{S_r} \otimes \rho_0^{S_{n-r}}$ where $\rho_0^{S_r}$ and $\rho_0^{S_{n-r}}$ are 
the trivial representations of $S_r$ and $S_{n-r}$, respectively.  

Due to the transitivity of induction,

\[ \rho_0 \uparrow_{S_r~\times~S_{n-r}}^{\mathbb{Z}_2~\wr~S_n} \ = \ \left\{ \rho_0 \uparrow_{S_r~\times~S_{n-r}}^{(\mathbb{Z}_2 
\wr S_r)~\times~(\mathbb{Z}_2~\wr~S_{n-r})} \right\} \uparrow_{(\mathbb{Z}_2~\wr~S_r)~\times~(\mathbb{Z}_2~\wr~S_{n-r})}^{\mathbb{Z}_2 
\wr S_n}. \]

\noindent
Thus since $\left\{ \rho_0^{S_r} \otimes \rho_0^{S_{n-r}} \right\} \uparrow_{S_r~\times~S_{n-r}}^{(\mathbb{Z}_2~\wr~S_r)~\times~
(\mathbb{Z}_2~\wr~S_{n-r})} \ = \ \left\{ \rho_0^{S_r} \uparrow_{S_r}^{\mathbb{Z}_2~\wr~S_r} \right\} \ \otimes \ \left\{ \rho_0^{S_{n-r}} 
\uparrow_{S_{n-r}}^{\mathbb{Z}_2~\wr~S_{n-r}} \right\}$, it follows that

\[ L(X) \ \ \cong \ \ \left\{ \rho_0^{S_r} \uparrow_{S_r}^{\mathbb{Z}_2~\wr~S_r} \ \otimes \ \rho_0^{S_{n-r}} 
\uparrow_{S_{n-r}}^{\mathbb{Z}_2~\wr~S_{n-r}} \right\} \uparrow_{(\mathbb{Z}_2~\wr~S_r)~\times~(\mathbb{Z}_2~\wr~S_{n-r})}^{\mathbb{Z}_2 
\wr S_n}. \]

It is a consequence of Corollary 4.4.7 of Greenhalgh (1989) that

\[ \rho_0^{S_r} \uparrow_{S_r}^{\mathbb{Z}_2~\wr~S_r} \ \ = \ \ \bigoplus_{i=0}^r \rho_{([i]; [r-i])} \ \ \ \mbox{and} \ \ \ 
\rho_0^{S_{n-r}} \uparrow_{S_{n-r}}^{\mathbb{Z}_2~\wr~S_{n-r}} \ \ = \ \ \bigoplus_{j=0}^{n-r} \rho_{([j]; [(n-r)-j])}, \]

\noindent
where $\rho_{([i]; [r-i])}$ is the irreducible representation of $\mathbb{Z}_2~\wr~S_r$ corresponding to the
trivial partitions $[i]$ and $[r-i]$ of $i$ and $r-i$, respectively, and $\rho_{([j]; [(n-r)-j])}$ is the irreducible
representation of $\mathbb{Z}_2\wr S_{n-r}$ corresponding to the trivial partitions $[j]$ and $[n-r-j]$ of $j$ and $(n-r)-j$,
respectively.

These results combine to show that

\[ \begin{array}{rcl}

L(X) &  \cong  &  \displaystyle \Bigg\{ \Big\{ \bigoplus_{i=0}^r \ \ \rho_{([i]; [r-i])} \Big\} \ \otimes \ \Big\{
\bigoplus_{j=0}^{n-r} \rho_{([j]; [(n-r)-j])} \Big\} \Bigg\} \uparrow_{(\mathbb{Z}_2~\wr~S_r)~\times~(\mathbb{Z}_2~\wr~
S_{n-r})}^{\mathbb{Z}_2~\wr~S_n} \vspace{1pc} \\

&  =  &  \displaystyle \bigoplus_{i=0}^r \ \bigoplus_{j=0}^{n-r} \ \left\{ \rho_{([i]; [r-i])} \ \otimes \ \rho_{([j];
[(n-r)-j])} \right\}\uparrow_{(\mathbb{Z}_2~\wr~S_r)~\times~(\mathbb{Z}_2~\wr~S_{n-r})}^{\mathbb{Z}_2~\wr~S_n}.

\end{array} \]

It follows from (the proof of) Lemma 4.4.5 of Greenhalgh (1989), which is a consequence of the ``inducing-up rule'' of Tokuyama
(1984), that 

\[ \left\{ \rho_{([i]; [r-i])} \ \otimes \ \rho_{([j]; [(n-r)-j])} \right\} \uparrow_{(\mathbb{Z}_2~\wr~S_r)~\times~(\mathbb{Z}_2 
\wr S_{n-r})}^{\mathbb{Z}_2~\wr~S_n} \ \ = \ \ \bigoplus_{[\lambda]} \ \bigoplus_{[\mu]} \ \rho_{([\lambda]; [\mu])}, \]

\noindent
where the range of summation over partitions $[\lambda]$ of $(i+j)$ is the range in the right-hand side of

\[ \left\{ \rho_{[i]} \otimes \rho_{[j]} \right\} \uparrow_{S_i~\times~S_j}^{S_{i+j}} \ = \ \bigoplus_{[\lambda]} \rho_{[\lambda]}
\]

\noindent
and, similarly, the range of summation over partitions $[\mu]$ of $n-(i+j)$ is the range in the right-hand side of

\[ \left\{ \rho_{[r-i]} \otimes \rho_{[(n-r)-j]} \right\} \uparrow_{S_{r-i}~\times~S_{(n-r)-j}}^{S_{n-(i+j)}} \ = \ 
\bigoplus_{[\mu]} \rho_{[\mu]}. \]

\noindent
It follows from Corollary 4.4.7 of Greenhalgh (1989) that 

\[ \begin{array}{c}

\displaystyle \left\{ \rho_{[i]} \otimes \rho_{[j]} \right\} \uparrow_{S_i~\times~S_j}^{S_{i+j}} \ \ = \ \ \bigoplus_{\ell=0}^{i
\wedge j} \rho_{[(i+j)-\ell,\ell]} \ \ \ \mbox{and} \vspace{1pc} \\

\displaystyle \left\{ \rho_{[r-i]} \otimes \rho_{[(n-r)-j]} \right\} \uparrow_{S_{r-i}~\times~S_{(n-r)-j}}^{S_{n-(i+j)}} \ \ = \ \
\bigoplus_{m=0}^{(r-i) \wedge ((n-r)-j)} \rho_{[n-(i+j)-m,m]}.

\end{array} \]

\noindent
These results combine to show that

\[ L(X) \ \ \cong \ \ \bigoplus_{i=0}^r \ \bigoplus_{j=0}^{n-r} \ \bigoplus_{\ell=0}^{i \wedge j} \ \bigoplus_{m=0}^{(r-i) \wedge
((n-r)-j)} \rho_{([(i+j)-\ell,\ell]; [n-(i+j)-m,m])}. \ \ \ \qed \]

In Sections~\ref{3.3} and \ref{3.4}, it will be more convenient to use the following decomposition of $L(X)$, which is a direct
consequence of combining $i$ and $j$ and changing the order of summation in Lemma~\ref{3.2.1}.

\begin{corollary} \label{3.2.2}
Let $L(X)$ be the set of all (complex-valued) functionals defined on $X = (\mathbb{Z}_2~\wr~S_n) /
(S_r~\times~S_{n-r})$.  Then, as a representation of $\mathbb{Z}_2~\wr~S_n$,

\[ L(X) \ \ \cong \ \ \bigoplus_{j=0}^n \ \bigoplus_{\ell=0}^{\lfloor j/2 \rfloor} \ \bigoplus_{i=\ell \vee (r-(n-j))}^{r \wedge
(j-\ell)} \ \bigoplus_{m=0}^{(r-i) \wedge ((n-j)-(r-i))} \rho_{([j-\ell,\ell]; [(n-j)-m,m])} \]

\noindent
where $\rho_{([j-\ell,\ell]; [(n-j)-m,m])}$ is the irreducible representation of $\mathbb{Z}_2~\wr~S_n$ corresponding to the
partitions $[j-\ell,\ell]$ and $[n-j-m,m]$ of $j$ and $(n-j)$, respectively.
\end{corollary}

The number of times that a particular representation $\rho$ occurs in the direct sums in Lemma~\ref{3.2.1} and Corollary~\ref{3.2.2}
is its multiplicity $m_{\rho}$ in the decomposition of $L(X)$.

\subsection{Analysis of the Independent Flips Model} \label{3.3}

In order to apply the Upper Bound Lemma (\ref{2.2.4}), we must now calculate the Fourier transform at each nontrivial irreducible
representation of $\mathbb{Z}_2~\wr~S_n$ occurring in the decomposition of $L(X)$.  We use the same technique as was used in Section
\ref{2.4}.

\begin{lemma} \label{3.3.1}
Let $P$ be the probability measure on $\mathbb{Z}_2~\wr~S_n$ defined in (\ref{3.1.1}).  Let $\rho = 
\rho_{([j-\ell,\ell]; [(n-j)-m,m])}$ be the nontrivial irreducible representation of $\mathbb{Z}_2~\wr~S_n$, corresponding to
the partitions $[j-\ell,\ell]$ and $[(n-j)-m,m]$ of $j$ and $n-j$, respectively, and occurring in the decomposition of $L(X)$.  Then, 
in a certain basis, the Fourier transform is

\[ \widehat{P}(\rho) \ \ = \ \ \left[ \frac{j^2}{n^2} \ - \ \frac{2\ell(j-\ell+1)}{n^2} \right] \widehat{I}(\rho) \]

\noindent
where $\widehat{I}(\rho)$ is the $d_{\rho}~\times~d_{\rho}$ matrix \onespace $\left[ \begin{array}{cc} I & 0 \\ 0 & 0 \\ \end{array}
\right]$ \twospace \textit{and $I$ is the $m_{\rho}$-dimensional identity matrix}.

\end{lemma}

\proof{Proof} Recall that the trivial representation of $\mathbb{Z}_2~\wr~S_n$ corresponds to the partition $[n]$ of $n$.  Thus we must
calculate the Fourier transform for the other nontrivial irreducible representations of $\mathbb{Z}_2~\wr~S_n$ occurring in the 
decomposition of $L(X)$, which were found in Corollary~\ref{3.2.2}. Notice that

\[ \begin{array}{rcl}

\widehat{P}(\rho)  &  =  &  \displaystyle \sum_{(\vec{0};\kappa) \in K} P(\vec{0};\kappa) \rho(\vec{0};\kappa) \ \ + \ \ 
\sum_{(\vec{u};\kappa) \in U_r K} P(\vec{u};\kappa) \rho(\vec{u};\kappa) \ \ + \ \ \sum_{(\vec{u};\kappa) \in U_{n-r} K}
P(\vec{u};\kappa) \rho(\vec{u};\kappa)  \vspace{1pc} \\

&  +  &  \displaystyle \sum_{(\vec{v};\tau\tau^{-1}\kappa) \in T_2 K \setminus U K} P(\vec{v};\tau\tau^{-1}\kappa)
\rho(\vec{v};\tau\tau^{-1}\kappa) \ \ + \ \ \sum_{(\vec{v};\tau \kappa) \in T_3 K} P(\vec{v};\tau \kappa) \rho(\vec{v};\tau \kappa).

\end{array} \]

Choose an orthonormal basis in $V$ such that the first $m_{\rho}$ basis vectors are $K$-invariant, as described in Section~\ref{2.2}.
It then follows from Lemma~\ref{2.2.3} that, in this basis,

\[ \begin{array}{rcl}

\displaystyle \sum_{(\vec{0};\kappa) \in K} P(\vec{0};\kappa) \rho(\vec{0};\kappa)  &  =  &  \displaystyle \frac{n(n+1) -
2r(n-r)}{4n^2 r!(n-r)!} \sum_{(\vec{0};\kappa) \in K} \rho(\vec{0};\kappa) \vspace{1pc} \\
&  =  &  \displaystyle \left[ \frac{n(n+1) - 2r(n-r)}{4n^2} \right] \cdot \left[ \frac{1}{r!(n-r)!} \right] \sum_{(\vec{0};\kappa)
\in K} \rho(\vec{0};\kappa) \vspace{1pc} \\
&  =  &  \displaystyle \left[ \frac{n(n+1) - 2r(n-r)}{4n^2} \right] \widehat{I}(\rho).

\end{array} \]

Recall that $U K \subseteq T_2 K$.  This is a result of the fact that for two-thirds of the elements in $T_2$, namely the odd
signed transpositions, $\vec{v} \in \mathbb{Z}_2^n$ has exactly one nonzero element.  Let $(\vec{u_i};e) \in U$ be the signed
identity whose only nonzero entry in $\vec{u_i}$ occurs in position $i$.  Let $(\vec{v_i};\tau_{ij}) \in T_2$ be the odd
transposition which transposes the entries in positions $i$ and $j$ and whose only nonzero entry in $\vec{v_i}$ occurs in position
$i$.

Notice that for any $\kappa \in K$, $(\vec{u_i};e) (\vec{0};\kappa) = (\vec{u_i}; \kappa) = 
(\vec{v_i};\tau_{ij}\tau_{ij}^{-1}\kappa) = (\vec{v_i};\tau_{ij}) (\vec{0};\tau_{ij}^{-1}\kappa)$ for $i \neq j$.  So it follows
that

\[ \begin{array}{l}

\displaystyle \sum_{(\vec{u};\kappa) \in U_r K} P(\vec{u};\kappa) \rho(\vec{u};\kappa) \ \ = \ \ \displaystyle \frac{r}{2n^2
r!(n-r)!} \sum_{(\vec{u};\kappa) \in U_r K} \rho(\vec{u};\kappa) \vspace{1pc} \\

\ \ \ = \ \ \displaystyle \frac{1}{2n^2 r!(n-r)!} \sum_{i=1}^r \sum_{(\vec{0};\kappa) \in K} \left[ \rho(\vec{u_i};\kappa) \ + \
\sum_{j \neq i} \rho(\vec{v_i};\tau_{ij}\tau_{ij}^{-1}\kappa) \right] \vspace{1pc} \\

\ \ \ = \ \ \displaystyle \frac{1}{2n^2} \sum_{i=1}^r \left[ \rho(\vec{u_i};e) \ + \ \sum_{j \neq i} 
\rho(\vec{v_i};\tau_{ij}) \right] \cdot \left[ \frac{1}{r!(n-r)!} \right] \sum_{(\vec{0};\kappa) \in K} \rho(\vec{0};\kappa).

\end{array} \]

\noindent
Similarly,

\[ \begin{array}{l}

\displaystyle \sum_{(\vec{u};\kappa) \in U_{n-r} K} P(\vec{u};\kappa) \rho(\vec{u};\kappa) \vspace{1pc} \\

\ \ \ = \ \ \displaystyle \frac{1}{2n^2} \sum_{i=r+1}^n \left[ \rho(\vec{u_i};e) \ + \ \sum_{j \neq i} \rho(\vec{v_i};\tau_{ij})
\right] \cdot \left[ \frac{1}{r!(n-r)!} \right] \sum_{(\vec{0};\kappa) \in K} \rho(\vec{0};\kappa). 

\end{array}\]

\noindent
Thus, since

\[ \begin{array}{l}

\displaystyle \sum_{(\vec{v};\tau\tau^{-1}\kappa) \in T_2 K \setminus U K} P(\vec{v};\tau\tau^{-1}\kappa)
\rho(\vec{v};\tau\tau^{-1}\kappa) \vspace{1pc} \\

\ \ \ = \ \ \displaystyle \frac{1}{2n^2} \left[ \sum_{(\vec{v};\tau) \in T_2^+} \rho(\vec{v};\tau) \right] \cdot \left[
\frac{1}{r!(n-r)!} \right] \sum_{(\vec{0};\kappa) \in K} \rho(\vec{0};\kappa),

\end{array} \]

\noindent
where $T_2^+$ is the set of all even transpositions in $T_2$, it follows that

\[ \begin{array}{l}

\displaystyle \sum_{(\vec{u};\kappa) \in U_r K} P(\vec{u};\kappa) \rho(\vec{u};\kappa) \ \ + \ \ \sum_{(\vec{u};\kappa)
\in U_{n-r} K} P(\vec{u};\kappa) \rho(\vec{u};\kappa)  \vspace{1pc} \\
\displaystyle \ \ \ \ \ \ \ + \ \ \sum_{(\vec{v};\tau\tau^{-1}\kappa) \in T_2 K \setminus U K} P(\vec{v};\tau\tau^{-1}\kappa)
\rho(\vec{v};\tau\tau^{-1}\kappa)  \vspace{1pc} \\
\displaystyle \ \ = \ \ \frac{1}{2n^2} \left[ \sum_{(\vec{u};e) \in U} \rho(\vec{u};e) \ + \ \sum_{(\vec{v};\tau) \in T_2}
\rho(\vec{v};\tau) \right] \cdot \left[ \frac{1}{r!(n-r)!} \right] \sum_{(\vec{0};\kappa) \in K} \rho(\vec{0};\kappa).

\end{array} \]

\noindent
Furthermore, we have

\[ \displaystyle \sum_{(\vec{v};\tau \kappa) \in T_3 K} P(\vec{v};\tau \kappa) \rho(\vec{v};\tau \kappa) \ \ = \ \ \frac{1}{2n^2}
\left[ \sum_{(\vec{v};\tau) \in T_3} \rho(\vec{v};\tau) \right] \cdot \left[ \frac{1}{r!(n-r)!} \right] \sum_{(\vec{0};\kappa) \in
K} \rho(\vec{0};\kappa). \]

\noindent
These results combine to show that

\[ \begin{array}{rcl}

\widehat{P}(\rho) \ &  =  &  \ \displaystyle \left[ \frac{n(n+1) - 2r(n-r)}{4n^2} \right] \widehat{I}(\rho)  \vspace{1pc}\\

&  +  &  \ \displaystyle \frac{1}{2n^2} \left[ \sum_{(\vec{u};e) \in U} \rho(\vec{u};e) \ + \ \sum_{(\vec{v};\tau) \in T} 
\rho(\vec{v};\tau) \ - \ \sum_{(\vec{v};\tau) \in T_1} \rho(\vec{v};\tau) \right] \widehat{I}(\rho).

\end{array} \]

It was shown in Section 3.3 of Schoolfield (1999) that $T$ splits into two conjugacy classes in $\mathbb{Z}_2~\wr~S_n$, namely, the
even transpositions $(\vec{v};\tau^+)$ (which change neither or both of the charges) and the odd transpositions $(\vec{v};\tau^-)$
(which change exactly one of the charges). It then follows from Lemma~\ref{2.4.1} that

\[ \sum_{(\vec{v};\tau) \in T} \rho(\vec{v};\tau) \ \ = \ \ \left[ \mbox{$\frac{1}{2}$} |T| \ r^+(\rho) \ + \ \mbox{$\frac{1}{2}$}
|T| \ r^-(\rho) \right] \ I \ \ = \ \ n(n-1) \ \left[ r^+(\rho) + r^-(\rho) \right] \ I, \]

\noindent
where $r^+(\rho) := \chi_{\rho}(\vec{v};\tau^+) / d_{\rho}$, $r^-(\rho) := \chi_{\rho}(\vec{v};\tau^-) / d_{\rho}$, and $I$ is the
$d_{\rho}$-dimensional identity matrix.  It follows from Lemma 3.5.1 of Schoolfield (1999) that

\[ \begin{array}{rcl}

\displaystyle r^+(\rho) \ & = & \ \displaystyle \frac{j(j-1) r(\lambda) \ \ + \ \ (n-j)(n-j-1) r(\mu)}{n(n-1)}, \vspace{1pc} \\  

\displaystyle r^-(\rho) \ & = & \ \displaystyle \frac{j(j-1) r(\lambda) \ \ - \ \ (n-j)(n-j-1) r(\mu)}{n(n-1)},  

\end{array} \]

\noindent
where $r(\lambda) := \chi_{[\lambda]}(\tau)/d_{[\lambda]}$, with $[\lambda] = [j-\ell,\ell]$ and $\tau \in S_j$, and $r(\mu) :=
\chi_{[\mu]}(\tau)/d_{[\mu]}$, with $[\mu] = [n-j-m,m]$ and $\tau \in S_{n-j}$.  From these results it follows that 

\[ \sum_{(\vec{v};\tau) \in T} \rho(\vec{v};\tau) \ \ = \ \ 2j(j-1) \ r(\lambda) \ I. \]

Since $L(X) \cong \rho_0 \uparrow_K^{\mathbb{Z}_2~\wr~S_n}$, where $\rho_0$ is the trivial representation of $K$, it follows that 
for any $\rho$ occurring in the decomposition of $L(X)$, $\rho \downarrow_K^{\mathbb{Z}_2~\wr~S_n}$ is the direct sum of 
$d_{\rho}$ copies of $\rho_0$.  Thus since $T_1 \subseteq K$, we have

\[ \sum_{(\vec{v};\tau) \in T_1} \rho(\vec{v};\tau) \ \ = \ \ |T_1| \ I \ \ = \ \ \left[ \mbox{$\frac{1}{2}$} n(n-1) \ - \ r(n-r)
\right] I. \]

It was shown in Section 3.3 of Schoolfield (1999) that $U$ is a conjugacy class in $\mathbb{Z}_2~\wr~S_n$.  It follows from
Lemma~\ref{2.4.1} that

\[ \sum_{(\vec{u};e) \in U} \rho(\vec{u};e) \ \ = \ \ |U| \left[ \chi_{\rho}(\vec{u};e) / d_{\rho} \right] \ I, \]

\noindent
where the right hand side is calculated at any $(\vec{u};e) \in U$.  It follows from Lemma 3.5.1 of Schoolfield (1999) that 
$\chi_{\rho}(\vec{u};e) / d_{\rho} \ = \ (2j-n) / n$.  From these results it follows that

\[ \sum_{(\vec{u};e) \in U} \rho(\vec{u};e) \ \ = \ \ \left( 2j - n \right) \ I. \]

\noindent
These results combine to show that

\[ \begin{array}{rcl}

\widehat{P}(\rho)  &  =  &  \displaystyle \left[ \frac{n(n+1) - 2r(n-r)}{4n^2} \right] \widehat{I}(\rho)  \vspace{1pc} \\
&     &  \displaystyle + \ \ \frac{1}{2n^2} \left[ (2j - n) \ + \ 2j(j-1) \ r(\lambda) \ - \ \mbox{$\frac{1}{2}$} n(n-1) \ + \
r(n-r) \right] \widehat{I}(\rho)  \vspace{1pc} \\
&  =  &  \displaystyle \left[ \frac{j}{n^2} \ + \ \frac{j(j-1)}{n^2} r(\lambda) \right] \widehat{I}(\rho).

\end{array} \]

\noindent
Recall that it follows from Lemma~\ref{2.4.2} that

\[ r(\lambda) \ \ = \ \ \displaystyle \frac{1}{j(j-1)} \left[ (j-\ell)(j-\ell-1) \ + \ \ell(\ell-3) \right]. \]

\noindent
Therefore, for the nontrivial irreducible representation of $\mathbb{Z}_2~\wr~S_n$, corresponding to the partitions 
$[j-\ell,\ell]$ and $[(n-j)-m,m]$ of $j$ and $n-j$, respectively, and occurring in the decomposition of $L(X)$,

\[ \begin{array}{rcl}

\widehat{P}(\rho)  &  =  &  \displaystyle \left[ \frac{j \ + (j-\ell)(j-\ell-1) \ + \ \ell(\ell-3)}{n^2} \right] \widehat{I}(\rho)
\vspace{1pc} \\
&  =  &  \displaystyle \left[ \frac{j^2}{n^2} \ - \ \frac{2\ell(j-\ell+1)}{n^2} \right] \widehat{I}(\rho). \ \ \ \ \ \qed

\end{array} \]

We have now established the results necessary to prove Theorem~\ref{3.1.3}.

\proof{Proof of Theorem \ref{3.1.3}} Notice that the probability measure $P$ defined in (\ref{3.1.1}) is clearly right $K$-invariant.
(In fact, it is bi-$K$-invariant.)  By applying the result from Lemma~\ref{3.3.1} to the Upper Bound Lemma (\ref{2.2.4}), we find that

\[ \displaystyle \| \widetilde{P^{*k}} - \widetilde{U} \|_{\mbox{\rm \scriptsize TV}}^2 \ \ \leq \ \ \displaystyle \mbox{$\frac{1}{4} 
\left[ 2^n \cdot {n \choose r} \right]$} \| \widetilde{P^{*k}} - \widetilde{U} \|_2^2
\ \ = \ \ \displaystyle \mbox{$\frac{1}{4}$} \sum_{\rho} d_{\rho} \ m_{\rho} \ \left[ \frac{j^2}{n^2} \ - \
\frac{2\ell(j-\ell+1)}{n^2} \right]^{2k}, \]

\noindent
where the sum is taken over all nontrivial irreducible representations $\rho = \rho_{([j-\ell,\ell]; [(n-j)-m,m])}$ of $\mathbb{Z}_2
\wr S_n$ occurring at least once in the decomposition of $L(X)$.  The factor $m_{\rho}$ comes from the trace of $\widehat{I}(\rho)$.  
It follows from Lemma 3.5.1 of Schoolfield (1999) that $d_{\rho} = {n \choose j} d_{[j-\ell,\ell]} \cdot d_{[(n-j)-m,m]}$.  It then 
follows from Lemma~\ref{2.3.1} that

\[ d_{\rho} \ \ = \ \ \mbox{${n \choose j} \left[ {j \choose \ell} - {j \choose \ell-1} \right] \cdot \left[ {n-j \choose m} - {n-j
\choose m-1} \right]$}. \]

\noindent
For notational purposes, define

\[ F(j,\ell) \ := \ \left[ \frac{j^2}{n^2} \ - \ \frac{2\ell(j-\ell+1)}{n^2} \right]. \]

Rather than explicitly calculate $m_{\rho}$, we will sum over the entire range of indices in the decomposition of $L(X)$ found in
Corollary~\ref{3.2.2}, thereby including each representation $\rho$ a total of $m_{\rho}$ times in the summation.  Thus we have
\num \begin{equation} \label{3.3.2}
\begin{array}{l}

\displaystyle \| \widetilde{P^{*k}} - \widetilde{U} \|_{\mbox{\rm \scriptsize TV}}^2 \ \ \leq \ \ \mbox{$\frac{1}{4} 
\left[ 2^n \cdot {n \choose r} \right]$} \| \widetilde{P^{*k}} - \widetilde{U} \|_2^2  \vspace{1pc} \\

\ \ \ \displaystyle = \ \ \mbox{$\frac{1}{4}$} \sum_{j=0}^n \ \sum_{\ell=0}^{\lfloor j/2 \rfloor} \ \sum_{i=\ell \vee
(r-(n-j))}^{r \wedge (j-\ell)} \ \sum_{m=0}^{(r-i) \wedge ((n-j)-(r-i))} \mbox{${n \choose j} \left[ {j \choose \ell} - {j
\choose \ell -1} \right] \cdot \left[ {n-j \choose m} - {n-j \choose m-1} \right] \cdot F(j,\ell)^{2k}$}  \vspace{1pc} \\

\ \ \ \displaystyle = \ \ \mbox{$\frac{1}{4}$} \sum_{j=0}^n \mbox{${n \choose j}$} \sum_{\ell=0}^{\lfloor j/2 \rfloor} \mbox{$\left[
{j \choose \ell} - {j \choose \ell -1} \right] \cdot F(j,\ell)^{2k}$} \ \sum_{i=\ell \vee (r-(n-j))}^{r \wedge (j-\ell)} \ 
\sum_{m=0}^{(r-i) \wedge ((n-j)-(r-i))} \mbox{$\left[ {n-j \choose m} - {n-j \choose m-1} \right]$}  \vspace{1pc} \\

\ \ \ \displaystyle = \ \ \mbox{$\frac{1}{4}$} \sum_{j=0}^n \mbox{${n \choose j}$} \sum_{\ell=0}^{\lfloor j/2 \rfloor} \mbox{$\left[
{j \choose \ell} - {j \choose \ell -1} \right] \cdot F(j,\ell)^{2k}$} \ \sum_{i=\ell \vee (r-(n-j))}^{r \wedge (j-\ell)} \
\mbox{${n-j \choose r-i}$}.

\end{array}
\end{equation}

\noindent
Notice that $\displaystyle \sum_{i=\ell \vee (r-(n-j))}^{r \wedge (j-\ell)} \ \mbox{${n-j \choose r-i}$} \ \leq \ \sum_{i=0}^{n-j}
\mbox{${n-j \choose i}$} \ = \ 2^{n-j}$. \  Thus we have

\[ \begin{array}{rcl}

\| \widetilde{P^{*k}} - \widetilde{U} \|_{\mbox{\rm \scriptsize TV}}^2  &  \leq  &  \mbox{$\frac{1}{4}$} \left[ 2^n \cdot {n \choose r}
\right] \| \widetilde{P^{*k}} - \widetilde{U} \|_2^2 \vspace{1pc} \\
&  \leq  &  \displaystyle \mbox{$\frac{1}{4}$} \sum_{j=0}^n 2^{n-j} {n \choose j} \sum_{\ell=0}^{\lfloor j/2 \rfloor}
\left[ {j \choose \ell} - {j \choose \ell -1} \right] \cdot \left[ \frac{j^2}{n^2} \ - \ \frac{2\ell(j-\ell+1)}{n^2} \right]^{2k}.

\end{array} \]

For each $1 \leq j \leq n$, it follows from (the calculations in the proof of) Theorem~\ref{2.5.3} that we may bound all but the $\ell
= 0$ term of the inner sum above by

\[ \begin{array}{l}

\displaystyle \sum_{\ell=1}^{\lfloor j/2 \rfloor} \left[ {j \choose \ell} - {j \choose \ell -1} \right] \cdot \left[ \frac{j^2}{n^2}
\ - \ \frac{2\ell(j-\ell+1)}{n^2} \right]^{2k}  \vspace{1pc} \\
\ \ \ = \ \displaystyle \left( \frac{j}{n} \right)^{4k} \sum_{\ell=1}^{\lfloor j/2 \rfloor} \left[ {j \choose \ell} - {j \choose
\ell-1} \right] \cdot \left[ 1 \ - \ \frac{2\ell(j-\ell+1)}{j^2} \right]^{2k} \ \ \leq \ \ \left( \frac{j}{n} \right)^{4k} 4 a^2 
e^{-c}

\end{array} \]

\noindent
for a universal constant $a > 0$, when $k \geq \frac{1}{4} j (\log j + c)$.  Since $n \geq j$, this is also true when $k \geq
\frac{1}{4} n (\log n + c)$.

We must also bound the term for the trivial representation $[\lambda] = [j]$ for $1 \leq j \leq n-1$.  Since in these cases $\ell =
0$, we have

\[ \left[ {j \choose \ell} - {j \choose \ell -1} \right] \cdot \left[ \frac{j^2}{n^2} \ - \ \frac{2\ell(j-\ell+1)}{n^2}
\right]^{2k} \ \ = \ \ \left( \frac{j}{n} \right)^{4k}. \]

\noindent
These results lead to the upper bound

\[ \begin{array}{rcl}

\| \widetilde{P^{*k}} - \widetilde{U} \|_{\mbox{\rm \scriptsize TV}}^2  &  \leq  &  \mbox{$\frac{1}{4}$} \left[ 2^n \cdot {n \choose r}
\right] \| \widetilde{P^{*k}} - \widetilde{U} \|_2^2 \vspace{1pc} \\
&  \leq  &  \displaystyle \mbox{$\frac{1}{4}$} \sum_{j=0}^n 2^{n-j} {n \choose j} \sum_{\ell=0}^{\lfloor j/2 \rfloor}
\left[ {j \choose \ell} - {j \choose \ell -1} \right] \cdot \left[ \frac{j^2}{n^2} \ - \ \frac{2\ell(j-\ell+1)}{n^2} \right]^{2k}
\vspace{1pc} \\
&  \leq  &  \displaystyle a^2 e^{-c} \sum_{j=1}^n 2^{n-j} {n \choose j} \left( \frac{j}{n} \right)^{4k} \ + \ \mbox{$\frac{1}{4}$}
\sum_{j=1}^{n-1} 2^{n-j} {n \choose j} \left( \frac{j}{n} \right)^{4k}.

\end{array} \]

\noindent
Now notice that, when $k = \frac{1}{4}n (\log n + c)$, then

\[ \displaystyle \left(\frac{j}{n}\right)^{4k} = \left(\frac{j}{n}\right)^{4\left[\frac{1}{4}n\log(n) + \frac{1}{4}cn\right]} =
\left(\frac{j}{n}\right)^{-n \left[-\log(n) - c\right]} = \left(\frac{e^{-c}}{n}\right)^{-n \log\left(\frac{j}{n}\right)}, \]

\noindent
which combines with the results above to give

\[ \begin{array}{rcl}

\| \widetilde{P^{*k}} - \widetilde{U} \|_{\mbox{\rm \scriptsize TV}}^2  &  \leq  &  \mbox{$\frac{1}{4}$} \left[ 2^n \cdot {n \choose r}
\right] \| \widetilde{P^{*k}} - \widetilde{U} \|_2^2 \vspace{1pc} \\
& \leq & \displaystyle a^2 e^{-c} \sum_{j=1}^n 2^{n-j} {n \choose j} \left(\frac{e^{-c}}{n}\right)^{-n
\log\left(\frac{j}{n}\right)} \vspace{1pc} \\
&   & \displaystyle + \ \ \mbox{$\frac{1}{4}$} \sum_{j=1}^{n-1} 2^{n-j} {n \choose j} \left(\frac{e^{-c}}{n}\right)^{-n 
\log\left(\frac{j}{n}\right)}.

\end{array} \]

\noindent
If we let $i = n-j$ it follows that

\[ \begin{array}{rcl}

\| \widetilde{P^{*k}} - \widetilde{U} \|_{\mbox{\rm \scriptsize TV}}^2  &  \leq  &  \mbox{$\frac{1}{4}$} \left[ 2^n \cdot {n \choose r}
\right] \| \widetilde{P^{*k}} - \widetilde{U} \|_2^2 \vspace{1pc} \\
& \leq & \displaystyle a^2 e^{-c} \sum_{i=0}^{n-1} \frac{1}{i!} (2n)^i \left(\frac{e^{-c}}{n}\right)^i \ \ \ + \ \ \
\mbox{$\frac{1}{4}$} \sum_{i=1}^{n-1} \frac{1}{i!} (2n)^i \left(\frac{e^{-c}}{n}\right)^i \vspace{1pc} \\
& = & \displaystyle a^2 e^{-c} \sum_{i=0}^{n-1} \frac{1}{i!} \left(2e^{-c}\right)^i \ \ \ + \ \ \ \mbox{$\frac{1}{2}$} e^{-c}
\sum_{i=0}^{n-2} \frac{1}{(i+1)!} \left(2e^{-c}\right)^i \vspace{1pc} \\
& \leq & \displaystyle a^2 e^{-c} \exp\left( 2e^{-c} \right) \ \ \ + \ \ \ \mbox{$\frac{1}{2}$} e^{-c} \exp\left( 2e^{-c}
\right).   

\end{array} \]

\noindent
Since $c > 0$, we have $\exp(2e^{-c}) < e^2$.  Therefore

\[ \| \widetilde{P^{*k}} - \widetilde{U} \|_{\mbox{\rm \scriptsize TV}}^2 \ \ \leq \ \ \mbox{$\frac{1}{4} \left[ 2^n \cdot {n \choose r}
\right]$} \| \widetilde{P^{*k}} - \widetilde{U} \|_2^2 \ \ \leq \ \ \left[ \left( a^2 + \mbox{$\frac{1}{2}$} \right)
e^2 \right] e^{-c}, \]

\noindent
from which the desired result follows.  $\qed$

Theorem~\ref{3.1.3} shows that $k = \frac{1}{4} n \left(\log n + c\right)$ steps are sufficient for the (normalized) $\ell^2$
distance, and hence the total variation distance, to become small.  A lower bound in the (normalized) $\ell^2$ metric can also be
derived by examining $2n \left( 1 - \frac{1}{n} \right)^{4k}$, which is the dominant contribution to the summation (\ref{3.3.2}) 
from the proof of Theorem~\ref{3.1.3}.  This term corresponds to the choice $j = n-1$ and $\ell = 0$.  Notice that $k = \frac{1}{4}
n \left(\log n - c\right)$ steps are necessary for just this term to become small.

Recall that Theorem~\ref{2.5.5} shows that, for values of $r$ not too small, $k =   \frac{1}{4} n \left(\log n - c\right)$
steps are necessary for the total variation distance to uniformity to become small in our variant of the Bernoulli--Laplace
diffusion model.  This is exactly the independent flips model, if the charges of the balls are ignored.  For such $r$ (in
particular, for the symmetric case $r = n/2$), Theorem~\ref{2.5.5} provides a matching lower bound on the distance to uniformity in
the total variation metric.  The upper bound in Theorem~\ref{3.1.3}, taken together with its matching lower bound, gives another
example of the ``cutoff phenomenon.''  

\subsection{Analysis of the Paired Flips Model.} \label{3.4}

At each step of the process introduced in Section~\ref{3.1}, the charges of the switched balls were changed independently.  Suppose
instead that, at each step, the charges of the switched balls are either both changed or both not changed.  It is this type of
process that we now examine.

We now describe a slight variant of the independent flips model introduced in Section~\ref{3.1}.  Suppose that the balls and racks
are as described in Section~\ref{3.1}.  At each step, independently choose two positions $p$ and $q$ uniformly from $\{1, 2, \ldots,
n\}$.  

If $p \neq q$, switch the balls in positions $p$ and $q$.  Then independently, with probability $1/2$, change the charge of the ball
moved to position $p$ and change the charge of the ball moved to position $q$.  Notice that this combination of operations is an
even transposition; the probability that an odd transposition occurs now vanishes.  Then, if necessary, permute the balls on each of
the two racks so that their labels are in increasing order.  

If $p = q$ (which occurs with probability $1/n$), leave the balls in their current positions.  Then, again independently with
probability $1/2$, change the charge of the ball in position $p = q$.  The probabilities of the identity and of the signed 
identities are thus unchanged from the independent flips model.  We refer to the process described above as the \emph{paired flips}
model.

The following analysis, while similar in format, is different in content from that found in Section~\ref{3.1}.  Let $T$ be the set
of all even transpositions in $\mathbb{Z}_2~\wr~S_n$.  Also, let $T_1$ be the set of all even transpositions in $K$, let $T_2$ be
the set of all even transpositions in $(\mathbb{Z}_2~\wr~K) \setminus K$, and let $T_3$ be the set of all even transpositions
in $(\mathbb{Z}_2~\wr~S_n) \setminus (\mathbb{Z}_2~\wr~K)$.  Thus $T = T_1 \cup T_2 \cup T_3$.  Notice that $\vec{v} =
\vec{0} \in \mathbb{Z}_2^n$ for any $(\vec{v};\tau) \in T_1$, that $\vec{v} \in \mathbb{Z}_2^n$ has exactly two
nonzero entries for any $(\vec{v};\tau) \in T_2$, and that $\vec{v} \in \mathbb{Z}_2^n$ has zero or two nonzero entries for
any $(\vec{v};\tau) \in T_3$.  Finally, let $U$ be the set of all signed identities in $\mathbb{Z}_2~\wr~S_n$.  Notice that $U
\subseteq \mathbb{Z}_2~\wr~K$.  Recall that for any $(\vec{u};e) \in U$, $\vec{u} \in \mathbb{Z}_2^n$ has exactly one
nonzero entry.

As with the classical model, before the two positions $p$ and $q$ have been chosen, the balls in the first rack may be permuted  
among themselves and the balls in the second rack may be permuted among themselves, without changing the state of the system
following the switch.  Thus, at each step, a random element of $K$ is effectively generated whenever the procedure described above
results in the identity or in an even transposition in $T_1$; this event occurs with probability 

\[ \frac{1}{2n} + \frac{\frac{1}{2} n(n-1) - r(n-r)}{n^2} = \frac{n^2 - 2r(n-r)}{2n^2}. \]

A similar analysis shows that the procedure effectively generates a random element of $U K$ with probability $\frac{1}{2n}$, a 
random element of $T_2 K$ with probability $\frac{n(n-1) - 2r(n-r)}{2n^2}$, and a random element of $T_3 K$ with probability
$\frac{2r(n-r)}{n^2}$.

Notice that each element of $T_3 K$ can be uniquely written as $(\vec{v};\tau\kappa)$, where $(\vec{v};\tau) \in T_3$ and 
$(\vec{0};\kappa) \in K$, each element of $T_2 K$ can be uniquely written as $(\vec{v};\tau\tau^{-1}\kappa)$, where $(\vec{v};\tau)
\in T_2$ and $(\vec{0};\tau^{-1}\kappa) \in K$, and each element of $U K$ can be uniquely written as $(\vec{u};\kappa)$, where
$(\vec{u};e) \in U$ and $(\vec{0};\kappa) \in K$.  Thus

\[ \begin{array}{rcl}

|K|  &  =  &  r!(n-r)!, \vspace{.5pc} \\
|U K|  &  =  &  n \, r!(n-r)!, \vspace{.5pc} \\
|T_2 K|  &  =  &  \left[ \mbox{$\frac{1}{2}$} n(n-1) - r(n-r) \right] r!(n-r)!, \mathrm{\ and} \vspace{.5pc} \\
|T_3 K|  &  =  &  2r(n-r) \cdot r!(n-r)!.

\end{array} \]

This paired flip variant of the signed Bernoulli--Laplace diffusion model may be modeled formally by a probability measure $Q$ on
the hyperoctahedral group $\mathbb{Z}_2~\wr~S_n$.  We may thus define the following probability measure on the set of all signed
permutations of $\mathbb{Z}_2~\wr~S_n$:
\num \begin{equation} \label{3.4.1}
\begin{array}{rcll}

Q(\vec{0};\kappa) & := & \frac{n^2- 2r(n-r)}{2n^2 r!(n-r)!} & \mbox{where $(\vec{0};\kappa)
\in K$}, \vspace{.5pc} \\

Q(\vec{u};\kappa) & := & \frac{1}{2n^2 r!(n-r)!} & \mbox{where $(\vec{u};\kappa) \in U K$}, \vspace{.5pc} \\

Q(\vec{v};\tau\tau^{-1}\kappa) & := & \frac{1}{n^2 r!(n-r)!} & \mbox{where $(\vec{v};\tau\tau^{-1}\kappa) \in T_2 K$}, 
\vspace{.5pc} \\

Q(\vec{v};\tau\kappa) & := & \frac{1}{n^2 r!(n-r)!} & \mbox{where $(\vec{v};\tau\kappa) \in T_3 K$, and}
\vspace{.5pc} \\

Q(\vec{x};\pi) & := & 0 & \mbox{otherwise}.

\end{array}
\end{equation}

In order to continue our analysis of the paired flips model, we must now calculate the Fourier transform of $Q$ at each nontrivial
irreducible representation of $\mathbb{Z}_2~\wr~S_n$ occurring in the decomposition of $L(X)$.  We use the same technique as was
used in Section~\ref{3.3}.

\begin{lemma} \label{3.4.2}
Let $Q$ be the probability measure on $\mathbb{Z}_2~\wr~S_n$ defined in (\ref{3.4.1}).  Let $\rho = 
\rho_{([j-\ell,\ell]; [(n-j)-m,m])}$ be the nontrivial irreducible representation of $\mathbb{Z}_2~\wr~S_n$, corresponding
the partitions $[j-\ell,\ell]$ and $[(n-j)-m,m]$ of $j$ and $n-j$, respectively, and occurring in the decomposition of $L(X)$.  
Then, in a certain basis, the Fourier transform is

\[ \widehat{Q}(\rho) \ \ = \ \ \mbox{$\left[ \frac{j^2}{n^2} \ - \ \frac{2\ell(j-\ell+1)}{n^2} \ + \ \frac{(n-j)^2}{n^2} \ - \ 
\frac{2m(n-j-m+1)}{n^2} \ - \ \frac{n-j}{n^2} \right]$} \widehat{I}(\rho) \]

\noindent
where $\widehat{I}(\rho)$ is the $d_{\rho}~\times~d_{\rho}$ matrix \onespace $\left[ \begin{array}{cc} I & 0 \\ 0 & 0 \\ \end{array}
\right]$ \twospace \textit{and $I$ is the $m_{\rho}$-dimensional identity matrix}.

\end{lemma}

\proof{Proof} Recall that the trivial representation of $\mathbb{Z}_2~\wr~S_n$ corresponds to the partition $[n]$ of $n$.  Thus we must
calculate the Fourier transform for the other nontrivial irreducible representations of $\mathbb{Z}_2~\wr~S_n$ occurring in the 
decomposition of $L(X)$, which were found in Corollary~\ref{3.2.2}. Notice that

\[ \begin{array}{rcl}

\widehat{Q}(\rho)  &  =  &  \displaystyle \sum_{(\vec{0};\kappa) \in K} Q(\vec{0};\kappa) \rho(\vec{0};\kappa) \ \ + \ \ 
\sum_{(\vec{u};\kappa) \in U K} Q(\vec{u};\kappa) \rho(\vec{u};\kappa) \vspace{1pc} \\

&  +  &  \displaystyle \sum_{(\vec{v};\tau\tau^{-1}\kappa) \in T_2 K} Q(\vec{v};\tau\tau^{-1}\kappa) 
\rho(\vec{v};\tau\tau^{-1}\kappa) \ \ + \ \ \sum_{(\vec{v};\tau \kappa) \in T_3 K} Q(\vec{v};\tau \kappa) \rho(\vec{v};\tau \kappa).

\end{array} \]

Choose an orthonormal basis in $V$ such that the first $m_{\rho}$ basis vectors are $K$-invariant, as described in Section~\ref{2.2}.
It then follows from Lemma~\ref{2.2.3} that, in this basis,

\[ \begin{array}{rcl}

\displaystyle \sum_{(\vec{0};\kappa) \in K} Q(\vec{0};\kappa) \rho(\vec{0};\kappa)  &  =  &  \displaystyle \frac{n^2 - 2r(n-r)}{2n^2
r!(n-r)!} \sum_{(\vec{0};\kappa) \in K} \rho(\vec{0};\kappa) \vspace{1pc} \\
&  =  &  \displaystyle \left[ \frac{n^2 - 2r(n-r)}{2n^2} \right] \cdot \left[ \frac{1}{r!(n-r)!} \right] \sum_{(\vec{0};\kappa) \in
K} \rho(\vec{0};\kappa) \vspace{1pc} \\
&  =  &  \displaystyle \left[ \frac{n^2 - 2r(n-r)}{2n^2} \right] \widehat{I}(\rho).

\end{array} \]

\noindent
Notice that

\[ \begin{array}{l}

\displaystyle \sum_{(\vec{u}; \kappa) \in U K} Q(\vec{u}; \kappa) \rho(\vec{u}; \kappa) \ \ =  \ \ \displaystyle \frac{1}{2n^2}
\left[ \sum_{(\vec{u};e) \in U} \rho(\vec{u};e) \right] \cdot \left[ \frac{1}{r!(n-r)!} \right] \sum_{(\vec{0};\kappa) \in K}
\rho(\vec{0};\kappa), \vspace{1pc} \\

\displaystyle \sum_{(\vec{v};\tau\tau^{-1}\kappa) \in T_2 K} Q(\vec{v};\tau\tau^{-1}\kappa) \rho(\vec{v};\tau\tau^{-1}\kappa) \ \ =
\ \ \displaystyle \frac{1}{n^2} \left[ \sum_{(\vec{v};\tau) \in T_2} \rho(\vec{v};\tau) \right] \cdot \left[ \frac{1}{r!(n-r)!}
\right] \sum_{(\vec{0};\kappa) \in K} \rho(\vec{0};\kappa), \vspace{1pc} \\

\mathrm{and\ \ } \mbox{$\displaystyle \sum_{(\vec{v};\tau \kappa) \in T_3 K} Q(\vec{v};\tau \kappa) \rho(\vec{v};\tau \kappa) \ \ = \
\ \displaystyle \frac{1}{n^2} \left[ \sum_{(\vec{v};\tau) \in T_3} \rho(\vec{v};\tau) \right] \cdot \left[ \frac{1}{r!(n-r)!} \right]
\sum_{(\vec{0};\kappa) \in K} \rho(\vec{0};\kappa)$}.

\end{array} \]

\noindent
These results combine to show that

\[ \begin{array}{rcl}

\widehat{Q}(\rho)  &  =  &  \displaystyle \left[ \frac{n^2 - 2r(n-r)}{2n^2} \right] \widehat{I}(\rho)  \vspace{1pc}\\

&  &  \displaystyle + \ \ \frac{1}{n^2} \left[ \mbox{$\frac{1}{2}$} \sum_{(\vec{u};e) \in U} \rho(\vec{u};e) \ + \
\sum_{(\vec{v};\tau) \in T} \rho(\vec{v};\tau) \ - \ \sum_{(\vec{v};\tau) \in T_1} \rho(\vec{v};\tau) \right] \widehat{I}(\rho).

\end{array} \]

Recall that the even transpositions $T$ form a conjugacy class in $\mathbb{Z}_2~\wr~S_n$.  It then follows from Lemma~\ref{2.4.1} that

\[ \sum_{(\vec{v};\tau) \in T} \rho(\vec{v};\tau) \ \ = \ \ |T| \ r^+(\rho) \ I \ \ = \ \ n(n-1) \ r^+(\rho) \ I, \]

\noindent
where $r^+(\rho) := \chi_{\rho}(\vec{v};\tau^+) / d_{\rho}$ for any even transposition $(\vec{v};\tau^+)$ and $I$ is the   
$d_{\rho}$-dimensional identity matrix.  Recall that $r^+(\rho)$ was determined in the proof of Lemma~\ref{3.3.1}.  From these results
it follows that 

\[ \sum_{(\vec{v};\tau) \in T} \rho(\vec{v};\tau) \ \ = \ \ \left[ j(j-1) \ r(\lambda) \ + \ (n-j)(n-j-1) \ r(\mu) \right] \ I, \]

\noindent
where $r(\lambda) := \chi_{[\lambda]}(\tau) / d_{[\lambda]}$, with $[\lambda] = [j-\ell,\ell]$ and $\tau \in S_j$, and $r(\mu)
:= \chi_{[\mu]}(\tau) / d_{[\mu]}$ with $[\mu] = [(n-j)-m,m]$ and $\tau \in S_{n-j}$.  As in the proof of Lemma~\ref{3.3.1},

\[ \begin{array}{c}

\displaystyle \sum_{(\vec{v};\tau) \in T_1} \rho(\vec{v};\tau) \ \ = \ \ |T_1| \ I \ \ = \ \ \left[ \mbox{$\frac{1}{2}$} n(n-1) \ -
\ r(n-r) \right] I \ \ \ \mathrm{and} \vspace{.5pc}\\

\displaystyle \sum_{(\vec{u};e) \in U} \rho(\vec{u};e) \ \ = \ \ \left( 2j - n \right) \ I.

\end{array} \]

\noindent
These results combine to show that

\[ \begin{array}{rcl}

\widehat{Q}(\rho)  &  =  &  \displaystyle \left[ \frac{n^2 - 2r(n-r)}{2n^2} \right] \widehat{I}(\rho)  \vspace{1pc} \\
&     &  + \ \ \left[ \frac{ \mbox{$\frac{1}{2}$} (2j - n) \ + \ j(j-1) \ r(\lambda) \ + \ (n-j)(n-j-1) \ r(\mu) \ - \
\mbox{$\frac{1}{2}$} n(n-1) \ + \ r(n-r)}{n^2} \right] \displaystyle \widehat{I}(\rho)  \vspace{1pc} \\
&  =  &  \displaystyle \left[ \frac{j}{n^2} \ + \ \frac{j(j-1)}{n^2} r(\lambda) \ + \ \frac{(n-j)(n-j-1)}{n^2} r(\mu) \right] 
\widehat{I}(\rho).

\end{array} \]

\noindent
Recall that it follows from Lemma~\ref{2.4.2} that

\[ r(\lambda) \ \ = \ \ \displaystyle \frac{1}{j(j-1)} \left[ (j-\ell)(j-\ell-1) \ + \ \ell(\ell-3) \right] \ \ \ \mbox{and} \]

\[ r(\mu) \ = \ \displaystyle \frac{1}{(n-j)(n-j-1)} \left[ (n-j-m)(n-j-m-1) \ + \ m(m-3) \right]. \]

\noindent
Therefore, for the nontrivial irreducible representation of $\mathbb{Z}_2~\wr~S_n$, corresponding the partitions 
$[j-\ell,\ell]$ and $[(n-j)-m,m]$ of $j$ and $n-j$, respectively, and occurring in the decomposition of $L(X)$,

\[ \begin{array}{rcl}

\widehat{Q}(\rho)  &  =  &  \left[ \frac{j \ + (j-\ell)(j-\ell-1) \ + \ \ell(\ell-3) \ + \ (n-j-m)(n-j-m-1) \ + \ m(m-3)}{n^2}
\right] \displaystyle \widehat{I}(\rho) \vspace{1pc} \\
&  =  &  \left[ \frac{j^2}{n^2} \ - \ \frac{2\ell(j-\ell+1)}{n^2} \ + \ \frac{(n-j)^2}{n^2} \ - \ \frac{2m(n-j-m+1)}{n^2} \ - \
\frac{n-j}{n^2} \right] \displaystyle \widehat{I}(\rho). \ \ \ \ \ \qed

\end{array} \]

The following result establishes an upper bound on both the total variation distance and the $\ell^2$ distance between 
$\widetilde{Q^{*k}}$ and $\widetilde{U}$, where $\widetilde{Q^{*k}}$ is the probability measure on the homogeneous space $X = 
(\mathbb{Z}_2~\wr~S_n) / (S_r~\times~S_{n-r})$, induced by the convolution $Q^{*k}$ of $Q$ with itself $k$ times, and 
$\widetilde{U}$ is the uniform probability measure on $X$.

\begin{theorem} \label{3.4.3}
Let $Q$ and $U$ be the probability measures on the hyperoctahedral group $\mathbb{Z}_2~\wr~S_n$ defined in
(\ref{3.4.1}) and (\ref{3.1.2}), respectively.  Let $\widetilde{Q^{*k}}$ be the probability measure on the homogeneous space $X =
(\mathbb{Z}_2~\wr~S_n) / (S_r~\times~S_{n-r})$ induced by $Q^{*k}$ and let $\widetilde{U}$ be the uniform probability
measure defined on $X$.  Let $k = \frac{1}{2} n (\log n + c)$.  Then there exists a universal constant $b > 0$ such that

\[ \| \widetilde{Q^{*k}} - \widetilde{U} \|_{\mbox{\rm \scriptsize TV}} \ \ \leq \ \ \mbox{$\frac{1}{2} \left[ 2^n \cdot {n \choose
r} \right]^{1/2}$} \| \widetilde{Q^{*k}} - \widetilde{U} \|_2 \ \ \leq \ \ be^{-c/2} \ \ \ \mathrm{for\ all\ }
\mbox{$c > 0$}.
\]

\end{theorem}

\proof{Proof} Notice that the probability measure $Q$ defined in (\ref{3.4.1}) is clearly right $K$-invariant.  (In fact, it is
bi-$K$-invariant.)  By applying the result from Lemma~\ref{3.4.2} to the Upper Bound Lemma (\ref{2.2.4}), we find that

\[ \begin{array}{l}

\displaystyle \| \widetilde{Q^{*k}} - \widetilde{U} \|_{\mbox{\rm \scriptsize TV}}^2 \ \ \leq \ \ \displaystyle \mbox{$\frac{1}{4} 
\left[ 2^n \cdot {n \choose r} \right]$} \| \widetilde{Q^{*k}} - \widetilde{U} \|_2^2  \vspace{1pc} \\

\ \ \ = \displaystyle \ \ \mbox{$\frac{1}{4}$} \sum_{\rho} d_{\rho} \ m_{\rho} \ \mbox{$\left[ \frac{j^2}{n^2} \ - \
\frac{2\ell(j-\ell+1)}{n^2} \ + \ \frac{(n-j)^2}{n^2} \ - \ \frac{2m(n-j-m+1)}{n^2} \ - \ \frac{n-j}{n^2} \right]^{2k},$}

\end{array} \]

\noindent
where the sum is taken over all nontrivial irreducible representations $\rho = \rho_{([j-\ell,\ell]; [(n-j)-m,m])}$ of $\mathbb{Z}_2
\wr S_n$ occurring at least once in the decomposition of $L(X)$.  Recall from the proof of Theorem~\ref{3.1.3} that

\[ d_{\rho} \ \ = \ \ \mbox{${n \choose j} \left[ {j \choose \ell} - {j \choose \ell-1} \right] \cdot \left[ {n-j \choose m} - {n-j
\choose m-1} \right]$}. \]

As in the proof of Theorem~\ref{3.1.3}, rather than explicitly calculate $m_{\rho}$, we will sum over the entire range of indices in
the decomposition of $L(X)$ found in Corollary~\ref{3.2.2}, thereby including each representation $\rho$ a total of $m_{\rho}$ times
in the summation.  Thus we have
\num \begin{equation} \label{3.4.4}
\begin{array}{l}

\displaystyle \| \widetilde{Q^{*k}} - \widetilde{U} \|_{\mbox{\rm \scriptsize TV}}^2 \ \ \leq \ \ \mbox{$\frac{1}{4} 
\left[ 2^n \cdot {n \choose r} \right]$} \| \widetilde{Q^{*k}} - \widetilde{U} \|_2^2  \vspace{1pc} \\

\ \ \ \displaystyle = \ \ \Bigg\{ \mbox{$\frac{1}{4}$} \sum_{j=0}^n \ \sum_{\ell=0}^{\lfloor j/2 \rfloor} \ \sum_{i=\ell \vee 
(r-(n-j))}^{r \wedge (j-\ell)} \ \sum_{m=0}^{(r-i) \wedge ((n-j)-(r-i))} \mbox{${n \choose j} \left[ {j \choose \ell} - {j
\choose \ell -1} \right] \cdot \left[ {n-j \choose m} - {n-j \choose m-1} \right]$} \vspace{1pc} \\

\ \ \ \displaystyle \ \ \ \ \ \times \ \ \mbox{$\left[ \frac{j^2}{n^2} \ - \ \frac{2\ell(j-\ell+1)}{n^2} \ + \ \frac{(n-j)^2}{n^2} \
- \ \frac{2m(n-j-m+1)}{n^2} \ - \ \frac{n-j}{n^2} \right]^{2k}$} \Bigg\}.

\end{array}
\end{equation}

\noindent
Notice that, when $1 \leq j \leq n-1$,

\[ \begin{array}{l}

\left[ \frac{j^2}{n^2} \ - \ \frac{2\ell(j-\ell+1)}{n^2} \ + \ \frac{(n-j)^2}{n^2} \ - \ \frac{2m(n-j-m+1)}{n^2} \ - \ 
\frac{n-j}{n^2} \right]^{2k} \vspace{1pc} \\

\ \ \ = \ \left\{ \left(\frac{j}{n}\right)^2 \left[ 1 - \frac{2\ell(j-\ell+1)}{j^2} \right] \ + \ \left(\frac{n-j}{n}\right)^2
\left[ 1 - \frac{2m(n-j-m+1)}{(n-j)^2} - \frac{1}{n-j} \right] \right\}^{2k} \vspace{1pc} \\

\ \ \ \leq \ \max \left\{ \left(\frac{j}{n}\right)^{2k} \left[ 1 -  \frac{2\ell(j-\ell+1)}{j^2} \right]^{2k}, \ 
\left(\frac{n-j}{n}\right)^{2k} \left[ 1 - \frac{2m(n-j-m+1)}{(n-j)^2} - \frac{1}{n-j} \right]^{2k} \right\} \vspace{1pc} \\

\ \ \ \leq \ \left(\frac{j}{n}\right)^{2k} \left[ 1 -  \frac{2\ell(j-\ell+1)}{j^2} \right]^{2k} \ + \
\left(\frac{n-j}{n}\right)^{2k} \left[ 1 - \frac{2m(n-j-m+1)}{(n-j)^2} - \frac{1}{n-j} \right]^{2k}, 

\end{array} \]

\noindent
where the first inequality is due to the fact that $\left( \alpha x + (1-\alpha) y \right)^{2k} \leq \max \left\{ x^{2k}, y^{2k}
\right\}$, for $0 \leq \alpha \leq 1$.  It can be determined that, for all possible choices of $j$ and $m$ with $1 \leq j \leq n-1$
and $m \leq (n-j)/2$,

\[ \mbox{$1 - \frac{2m(n-j-m+1)}{(n-j)^2} - \frac{1}{n-j}$} \ \geq \ 0, \]

\noindent
except when $j = n-2$ and $m=1$.  But when $j = n-2$, notice that

\[ \sum_{m=0}^1 \mbox{$\left[ 1 - \frac{2m(n-j-m+1)}{(n-j)^2} - \frac{1}{n-j} \right]^{2k}$} \ \leq \ \sum_{m=0}^1 \mbox{$\left[ 1 -
\frac{2m(n-j-m+1)}{(n-j)^2} \right]^{2k}$}. \]

\noindent
These results combine to give the upper bound

\[ \begin{array}{l}

\displaystyle \| \widetilde{Q^{*k}} - \widetilde{U} \|_{\mbox{\rm \scriptsize TV}}^2 \ \ \leq \ \ \displaystyle \mbox{$\frac{1}{4} 
\left[ 2^n \cdot {n \choose r} \right]$} \| \widetilde{Q^{*k}} - \widetilde{U} \|_2^2  \vspace{1pc} \\

\ \ \ \displaystyle \leq \ \ \Bigg\{ \mbox{$\frac{1}{4}$} \sum_{j=1}^n \ \sum_{\ell=0}^{\lfloor j/2 \rfloor} \ \sum_{i=\ell \vee 
(r-(n-j))}^{r \wedge (j-\ell)} \ \sum_{m=0}^{(r-i) \wedge ((n-j)-(r-i))} \mbox{${n \choose j} \left[ {j \choose \ell} - {j \choose
\ell -1} \right] \cdot \left[ {n-j \choose m} - {n-j \choose m-1} \right]$} \vspace{1pc} \\

\ \ \ \displaystyle \ \ \ \ \ \times \ \ \mbox{$\left(\frac{j}{n}\right)^{2k} \left[ 1 -  \frac{2\ell(j-\ell+1)}{j^2} \right]^{2k}$}
\Bigg\} \vspace{1pc} \\

\ \ \ \displaystyle + \ \ \Bigg\{ \mbox{$\frac{1}{4}$} \sum_{j=0}^{n-1} \ \sum_{\ell=0}^{\lfloor j/2 \rfloor} \ \sum_{i=\ell \vee 
(r-(n-j))}^{r \wedge (j-\ell)} \ \sum_{m=0}^{(r-i) \wedge ((n-j)-(r-i))} \mbox{${n \choose j} \left[ {j \choose \ell} - {j \choose
\ell -1} \right] \cdot \left[ {n-j \choose m} - {n-j \choose m-1} \right]$} \vspace{1pc} \\

\ \ \ \displaystyle \ \ \ \ \ \times \ \ \mbox{$\left(\frac{n-j}{n}\right)^{2k} \left[ 1 - \frac{2m(n-j-m+1)}{(n-j)^2}
\right]^{2k}$} \Bigg\} \ \ + \ \ \mbox{$\frac{1}{4} \left( \frac{n-1}{n} \right)^{2k}$},

\end{array} \]

\noindent
where we must modify the first iterated summation to exclude the term for $j=n$ and $\ell=0$ (and hence $i=r$ and $m=0$) and we must
also modify the second iterated summation to exclude the term for $j=0$ and $m=0$ (and hence $\ell=0$ and $i=0$).  The final 
expression reintroduces the appropriate term in the exact formula for the squared $\ell^2$ distance at the representation
corresponding to $j=0$ and $m=0$.

Notice that if, in the second set of braces $\Bigg\{ \ \Bigg\}$, we put $i' = r-i$ and $j' = (n-r)-j$ and change the order of 
summation, we obtain

\[ \begin{array}{l}

\displaystyle \| \widetilde{Q^{*k}} - \widetilde{U} \|_{\mbox{\rm \scriptsize TV}}^2 \ \ \leq \ \ \displaystyle \mbox{$\frac{1}{4} 
\left[ 2^n \cdot {n \choose r} \right]$} \| \widetilde{Q^{*k}} - \widetilde{U} \|_2^2  \vspace{1pc} \\

\ \ \ \displaystyle \leq \ \ \Bigg\{ \mbox{$\frac{1}{4}$} \sum_{j=1}^n \ \sum_{\ell=0}^{\lfloor j/2 \rfloor} \ \sum_{i=\ell \vee 
(r-(n-j))}^{r \wedge (j-\ell)} \ \sum_{m=0}^{(r-i) \wedge ((n-j)-(r-i))} \mbox{${n \choose j} \left[ {j \choose \ell} - {j
\choose \ell -1} \right] \cdot \left[ {n-j \choose m} - {n-j \choose m-1} \right]$} \vspace{1pc} \\

\ \ \ \displaystyle \ \ \ \ \ \times \ \ \mbox{$\left(\frac{j}{n}\right)^{2k} \left[ 1 -  \frac{2\ell(j-\ell+1)}{j^2} \right]^{2k}$}
\Bigg\} \vspace{1pc} \\

\ \ \ \displaystyle + \ \ \Bigg\{ \mbox{$\frac{1}{4}$} \sum_{j=0}^{n-1} \ \sum_{m=0}^{\lfloor (n-j)/2 \rfloor} \ \sum_{i=m \vee 
(r-j)}^{r \wedge ((n-j)-m)} \ \sum_{\ell=0}^{(r-i) \wedge (j-(r-i))} \mbox{${n \choose j} \left[ {j \choose \ell} - {j
\choose \ell -1} \right] \cdot \left[ {n-j \choose m} - {n-j \choose m-1} \right]$} \vspace{1pc} \\

\ \ \ \displaystyle \ \ \ \ \ \times \ \ \mbox{$\left(\frac{n-j}{n}\right)^{2k} \left[ 1 - \frac{2m(n-j-m+1)}{(n-j)^2}
\right]^{2k}$} \Bigg\} \ \ + \ \ \mbox{$\frac{1}{4} \left( \frac{n-1}{n} \right)^{2k}$}.

\end{array} \]

\noindent
Notice that if we now put $j' = n-j$ and interchange the roles of $\ell$ and $m$, then the second summation becomes identical to the
first.  Thus, combining these summations and continuing as in the proof of Theorem~\ref{3.1.3}, we have

\[ \begin{array}{l}

\displaystyle \| \widetilde{Q^{*k}} - \widetilde{U} \|_{\mbox{\rm \scriptsize TV}}^2 \ \ \leq \ \ \displaystyle \mbox{$\frac{1}{4} 
\left[ 2^n \cdot {n \choose r} \right]$} \| \widetilde{Q^{*k}} - \widetilde{U} \|_2^2  \vspace{1pc} \\

\ \ \ \displaystyle \leq \ \ \Bigg\{ \mbox{$\frac{1}{2}$} \sum_{j=1}^n \ \sum_{\ell=0}^{\lfloor j/2 \rfloor} \ \sum_{i=\ell \vee 
(r-(n-j))}^{r \wedge (j-\ell)} \ \sum_{m=0}^{(r-i) \wedge ((n-j)-(r-i))} \mbox{${n \choose j} \left[ {j \choose \ell} - {j
\choose \ell -1} \right] \cdot \left[ {n-j \choose m} - {n-j \choose m-1} \right]$} \vspace{1pc} \\

\ \ \ \displaystyle \ \ \ \ \ \times \ \ \mbox{$\left(\frac{j}{n}\right)^{2k} \left[ 1 -  \frac{2\ell(j-\ell+1)}{j^2} \right]^{2k}$}
\Bigg\} \ \ + \ \ \mbox{$\frac{1}{4} \left( \frac{n-1}{n} \right)^{2k}$} \vspace{1pc} \\

\ \ \ \displaystyle \leq \ \ \mbox{$\frac{1}{2}$} \sum_{j=1}^n 2^{n-j} {n \choose j} \left(\frac{j}{n}\right)^{2k} 
\sum_{\ell=0}^{\lfloor j/2 \rfloor} \left[ {j \choose \ell} - {j \choose \ell -1} \right] \cdot \left[ 1 -
\frac{2\ell(j-\ell+1)}{j^2} \right]^{2k} \vspace{1pc} \\

\ \ \ \displaystyle \ \ \ \ \ + \ \ \mbox{$\frac{1}{4}$} \left( \frac{n-1}{n} \right)^{2k}.

\end{array} \]

Recall from the proof of Theorem~\ref{3.1.3} that, when $k \geq \frac{1}{4}n (\log n + c)$, we may bound the inner sum above using

\[ \left( \frac{j}{n} \right)^{2k} \sum_{\ell=0}^{\lfloor j/2 \rfloor} \left[ {j \choose \ell} - {j \choose \ell -1} \right] \cdot
\left[ 1 \ - \ \frac{2\ell(j-\ell+1)}{j^2} \right]^{2k} \ \ \leq \ \ \left( \frac{j}{n} \right)^{2k} 4 a^2 e^{-c} \ + \ \left(
\frac{j}{n} \right)^{2k} \]

\noindent
for $1 \leq j \leq n-1$, and using

\[ \left( \frac{j}{n} \right)^{2k} \sum_{\ell=1}^{\lfloor j/2 \rfloor} \left[ {j \choose \ell} - {j \choose \ell -1} \right] \cdot 
\left[ 1 \ - \ \frac{2\ell(j-\ell+1)}{j^2} \right]^{2k} \ \ \leq \ \ 4 a^2 e^{-c} \]

\noindent
for $j = n$.  So this is also true when $k = \frac{1}{2}n (\log n + c)$.  

Now notice that when $k = \frac{1}{2}n (\log n + c)$

\[ \left(\frac{j}{n}\right)^{2k} = \left( \frac{e^{-c}}{n} \right)^{-n \log(j/n)}. \]

\noindent
These results lead to the upper bound

\[ \begin{array}{rcl}

\| \widetilde{Q^{*k}} - \widetilde{U} \|_{\mbox{\rm \scriptsize TV}}^2  &  \leq  &  \mbox{$\frac{1}{4}$} \left[ 2^n \cdot {n \choose r}
\right] \| \widetilde{Q^{*k}} - \widetilde{U} \|_2^2 \vspace{1pc} \\

& \leq & \displaystyle 2a^2 e^{-c} \sum_{j=1}^n 2^{n-j} {n \choose j} \left(\frac{e^{-c}}{n}\right)^{-n
\log\left(\frac{j}{n}\right)} \vspace{1pc} \\

&   & \displaystyle + \ \ \mbox{$\frac{1}{2}$} \sum_{j=1}^{n-1} 2^{n-j} {n \choose j} \left(\frac{e^{-c}}{n}\right)^{-n 
\log\left(\frac{j}{n}\right)} \ \ + \ \ \mbox{$\frac{1}{4}$} \left( \frac{e^{-c}}{n} \right)^{-n \log \left(1 - \frac{1}{n}\right)}.

\end{array} \]

\noindent
Continuing as in the proof of Theorem~\ref{3.1.3}, we find that

\[ \begin{array}{rcl}

\| \widetilde{Q^{*k}} - \widetilde{U} \|_{\mbox{\rm \scriptsize TV}}^2  &  \leq  &  \mbox{$\frac{1}{4}$} \left[ 2^n \cdot {n \choose r}
\right] \| \widetilde{Q^{*k}} - \widetilde{U} \|_2^2 \vspace{1pc} \\

& \leq & \displaystyle 2a^2 e^{-c} \exp\left( 2e^{-c} \right) \ \ \ + \ \ \ e^{-c} \exp\left( 2e^{-c} \right) \ \ \ + \ \ \
\mbox{$\frac{1}{4}$} e^{-c}.

\end{array} \]

\noindent
Since $c > 0$, we have $\exp(2e^{-c}) < e^2$.  Therefore

\[ \| \widetilde{Q^{*k}} - \widetilde{U} \|_{\mbox{\rm \scriptsize TV}}^2 \ \ \leq \ \ \mbox{$\frac{1}{4} \left[ 2^n \cdot {n \choose r}
\right]$} \| \widetilde{Q^{*k}} - \widetilde{U} \|_2^2 \ \ \leq \ \ \left[ \left( 2a^2 + 1 \right) e^2 +
\mbox{$\frac{1}{4}$} \right] e^{-c}, \]

\noindent
from which the desired result follows.  $\qed$

Theorem~\ref{3.4.3} shows that $k = \frac{1}{2} n \left(\log n + c\right)$ steps are sufficient for the (normalized) $\ell^2$
distance, and hence the total variation distance, to become small.  A lower bound in the (normalized) $\ell^2$ metric can also be
derived by examining $2n \left( 1 - \frac{1}{n} \right)^{4k}$, which is the dominant contribution to the summation (\ref{3.4.4}) 
from the proof of Theorem~\ref{3.4.3}.  This term corresponds to the choice $j = n-1$ and $\ell = 0$.  Notice that $k = \frac{1}{4}
n \left(\log n - c\right)$ steps are necessary for just this term to become small.  Furthermore, our upper ($\frac{1}{2} n \log n$)
and lower ($\frac{1}{4} n \log n$) bounds on the number of steps required for the (normalized) $\ell^2$ distance to become small
differ by a constant factor.  We have not been able to close this gap.  

Recall that Theorem~\ref{2.5.5} shows that, for values of $r$ not too small, $k =  \frac{1}{4} n \left(\log n - c\right)$
steps are necessary for the total variation distance to uniformity to become small in our variant of the Bernoulli--Laplace 
diffusion model.  This is exactly the paired flips model, if the charges of the balls are ignored.  For such $r$ (in particular, for
the symmetric case $r = n/2$), Theorem~\ref{2.5.5} provides a lower bound (differing from the upper bound by only a constant factor)
on the distance to uniformity in the total variation metric, just as in Section~\ref{3.3}.

\section*{Acknowledgments.} 
This paper formed a portion of the author's Ph.D. dissertation in the Department of Mathematical Sciences at the Johns Hopkins
University.  The author wishes to thank his advisor Jim Fill, whose assistance was invaluable, particularly in the proof of
Theorem~\ref{2.5.5}.

\twospace

\Line{\hfill
\AOPaddress{Clyde H. Schoolfield, Jr.\\
Department of Statistics\\
Harvard University\\
One Oxford Street\\
Cambridge, Massachusetts 02138\\
e-mail:  {\tt clyde@stat.harvard.edu}}}

\end{document}